\documentclass[3p,,preprint,12pt]{elsarticle}
\usepackage{amsmath}
\usepackage{graphicx}
\usepackage{wrapfig}
\usepackage{epsfig}
\usepackage{epstopdf}
\usepackage{amsmath}
\usepackage{times}
\usepackage{epsfig}
\usepackage{amssymb}
\usepackage{array}
\usepackage{latexsym}
\usepackage{amsthm}
\usepackage{amsfonts}
\usepackage{nomencl}
\usepackage{lineno}

\makenomenclature
\usepackage{hyperref}
\usepackage{array,tabularx}
\usepackage{scalerel,stackengine}
\usepackage{graphics}
\usepackage{epstopdf}
\usepackage[english]{babel}

\DeclareGraphicsExtensions{.png,.pdf}
\stackMath
\numberwithin{equation}{section}
\date{}
\newtheorem{t1}{Theorem}[section]
\newtheorem{p1}{Proposition}[section]
\newtheorem{l1}{Lemma}[section]
\newtheorem{c1}{Corollary}[section]

\usepackage[utf8]{inputenc}
\begin{document}
\begin{frontmatter} 

\title{On the existence and regularity of solutions of semi-hyperbolic patches to 2-D Euler equations with van der Waals gas\tnoteref{mytitlenote}} 
\author{Rahul Barthwal}
\author{T. Raja Sekhar}
\address{Department of Mathematics, Indian Institute of Technology Kharagpur, Kharagpur,  India} 


\cortext[mycorrespondingauthor]{Corresponding author} 
\ead{trajasekhar@maths.iitkgp.ac.in} 


\begin{abstract}
This article is concerned in establishing the existence and regularity of solution of semi-hyperbolic patch problem for two-dimensional isentropic Euler equations with van der Waals gas. This type of solution appears in the transonic flow over an airfoil and Guderley reflection and is very common in the numerical solution of Riemann problems. We use the idea of characteristic decomposition and bootstrap method to prove the existence of global smooth solution which is uniformly $C^{1, \frac{1}{2}}$ continuous up to the sonic curve. We also prove that the sonic curve is $C^{1, \frac{1}{2}}$ continuous. Further, we show the formation of shock as an envelope for positive characteristics before reaching their sonic points. 
\end{abstract} 


\end{frontmatter}

\section{Introduction}
Cauchy problem in several space dimension for hyperbolic system of conservation laws is a very important but complicated open problem. A particular kind of Cauchy problem in two-dimensional case is the two-dimensional Riemann problem which consists of initial data that are constant along any ray passing through origin. The study of two-dimensional Riemann problems are very interesting and challenging in the context of two-dimensional hyperbolic system of conservation laws. A significant research has been done for the two-dimensional compressible Euler system and various other important models for a typical case of two-dimensional Riemann problems which is known as the four wave Riemann problem. The four wave Riemann problem is an initial value problem where the initial data are constant in each of the four quadrants of the physical plane. For the two-dimensional compressible Euler system, a beautiful conjecture for the possible structures of solution for four wave Riemann problem was provided in the ground breaking paper of Zhang and Zheng \cite{zhang1989conjecture}. Several numerical results have been obtained in the field of two-dimensional Riemann problems for gasdynamics equations and many small-scale structures have been observed in those numerical simulations \cite{glimm2008transonic, schulz1993numerical, kurganov2002solution}. One of such structures is the semi-hyperbolic patch which appears very often in many cases of two-dimensional Riemann problems for compressible Euler system, pressure-gradient equations, magnetohydrodynamics and etc(For details see \cite{zheng2012systems, li1998two, song2009semi, li2011semi, chen2019semi}). The semi-hyperbolic patch is defined as a patch kind of solution in which one set of characteristics starts on a sonic curve and ends on either a transonic shock wave or a sonic curve. These type of solutions appear in many other situations too such as reflection of rarefaction wave along a compressive corner \cite{sheng2010critical}, transonic flow over an airfoil \cite{courant0} and Guderley shock reflection of the von Neumann triple point paradox \cite{tesdall2008self, tesdall2007triple}. These patch type solutions are very meaningful and important for the construction of global solution of mixed-type equations in future.

The semi-hyperbolic patch, first time, has been identified among the small-scale structures in the work of Song and Zheng \cite{song2009semi} for pressure-gradient system. The same problem for isothermal Euler equations was studied by Hu et al. in \cite{hu2012degenerate} while for isentropic case by Li and Zheng in \cite{li2011semi} and extended to magnetohydrodynamic system by Chen and Lai in \cite{chen2019semi}. The regularity of the solution of semi-hyperbolic patch problems have also been widely discussed. The regularity of solution of the semi-hyperbolic patch problem for the pressure-gradient system was discussed in \cite{wang2014regularity} and for isentropic Euler system in \cite{song2015regularity}. The regularity results for isothermal Euler equations have been discussed by Hu et al. \cite{hu2019regularity}. The regularity results obtained for isentropic Euler system in \cite{song2015regularity} were improved by Hu et al. in \cite{hu2018improved}. For more details on ongoing research on two-dimensional Euler system and related models, we refer the reader to \cite{li2010interaction, li2009interaction, hu2020sonic, lai2015centered, hu2014semi, li2006simple, li2011characteristic, zhang2014sonic, zhang2017existence, zhang2016structure, hu2020global, li2019degenerate, hu2019sonic, song2020regularity, lai2016expansion, sheng2018interaction}. All the above works are based on the beautiful concept of characteristic decomposition initiated in the work of Dai and Zhang \cite{dai2000existence}. The progress made in the field of semi-hyperbolic patch problems leads to a natural question of determining whether these results can be extended for more realistic gases, for instance, van der Waals gas. The main purpose of this article is to establish the existence and regularity results for the solution of the semi-hyperbolic patch problem for two-dimensional isentropic Euler equations with van der Waals gas. 
 
Let us consider the two-dimensional isentropic compressible Euler equations \cite{song2015regularity} as follows:
\begin{equation}\label{eq: 1}
\begin{aligned}
            \rho_t+(\rho u)_x+(\rho v)_y=0,\\
            (\rho u)_t+(\rho u^2+p)_x+(\rho u v)_y=0,\\
            (\rho v)_t+(\rho u v)_x+(\rho v^2+p)_y=0,
\end{aligned}
\end{equation}
where $\rho$ denotes the density, $u$ and $v$ denotes the flow velocity in the $x$ and $y$ direction, respectively and $p$ denotes the pressure of the gas.

We consider a polytropic van der Waals gas with the equation of state \cite{callen1998thermodynamics} as $p(\tau)= \frac{K}{(\tau-b)^{\gamma+1}}-\frac{a}{\tau^2}$ where $\gamma$ is a constant such that $\gamma \in (0, 1)$. The quantity $\tau=\frac{1}{\rho}$ is known as the specific volume of the gas, $K$ is a positive constant which depends on the entropy of the system. The attraction between the gas molecules is represented by the positive constant $a$ and the compressibility limit of these molecules in the gas is represented by the positive constant $b$. For the case $a=0$, this corresponds to dusty gas while for $a=b=0$, this behaves as polytropic ideal gas.

 The expression for speed of sound is given by $c(\tau)=\sqrt{-\tau^2 p'(\tau)}$ with 
 \begin{align}\label{eq: 2}
   p'(\tau)=-\dfrac{K(\gamma+1)}{(\tau-b)^{\gamma+2}}+\dfrac{2a}{\tau^3},~~~~~~~~~~~~~~~p''(\tau)=\dfrac{K(\gamma+1)(\gamma+2)}{(\tau-b)^{\gamma+3}}-\dfrac{6a}{\tau^4}.  
 \end{align}
 We give the following list of notations, these are very important in the further discussion of the paper
 \begin{align}\label{eq: 4}
 \begin{cases}
  \kappa(\tau)=\dfrac{-2p'(\tau)}{2p'(\tau)+\tau p''(\tau)}=\dfrac{2-\frac{2b}{\tau}-\frac{4a(\tau-b)^{\gamma+3}}{K(\gamma+1)\tau^4}}{\gamma+\frac{2b}{\tau}+\frac{2a(\tau-b)^{\gamma+3}}{K(\gamma+1)\tau^4}}, \vspace{0.15 cm} \\
  m(\tau)=\dfrac{\kappa(\tau)-1}{\kappa(\tau)+1}=\dfrac{2-\gamma-\frac{4b}{\tau}-\frac{2a(\tau-b)^{\gamma+3}}{K(\gamma+1)\tau^4}}{\gamma+2+\frac{6a(\tau-b)^{\gamma+3}}{K(\gamma+1)\tau^4}},  \vspace{0.15 cm} \\
    \mu^2(\tau)=\dfrac{1}{1+\kappa(\tau)}, ~~~\Omega(\tau, \omega)=m(\tau)-\tan^2 \omega, 
 \vspace{0.15 cm} \\
     \bar{\partial}_+= \cos \alpha \partial_\xi+\sin \alpha \partial_\eta, ~~
     \bar{\partial}_-= \cos \beta \partial_\xi+\sin \beta \partial_\eta, \vspace{0.15 cm} \\
     \sigma=\dfrac{\alpha+\beta}{2},~~~~\omega=\dfrac{\alpha-\beta}{2},~~~\tan \alpha=\lambda_+,~~\tan \beta=\lambda_- 
     \end{cases}
 \end{align}
 where $\alpha$ and $\beta$ are defined as characteristic angles and $\bar{\partial}_\pm$ denotes the normalized directional derivatives along the characteristic directions in self-similar plane \cite{li2011characteristic}.
 
 From the expressions above we can observe that for sufficiently large $\tau_1>b$, when $\tau>\tau_1$ we have 
 $p'(\tau)<0,~ p''(\tau)>0,~ \kappa(\tau)>0$ and $0<m(\tau)<1$. Further, for some technical reasons we adopt the hypothesis that $\kappa'(\tau)>0$ without loss of generality.
 
 This paper is organized as follows. In Section $\ref{2}$ we give some preliminaries and mainly interested with characteristic decompositions in terms of characteristic angles and the speed of sound. We define our problem precisely and establish the boundary data estimates to prove the existence of local solution in Section $\ref{3}$. Section $\ref{4}$ is devoted to construct the uniform lower and upper bounds of the characteristic directional derivatives of speed of sound. We discuss the global existence of solution by extending the local solution up to the sonic boundary by solving many small Goursat problems in each step of extension in Section $\ref{5}$. In Section $\ref{6}$ we study the formation of shock as an envelope for positive characteristics before reaching their sonic points. The regularity of solution in partial hodograph plane and self-similar plane is established in Section $\ref{7}$ and $\ref{8}$, respectively. In Section $\ref{9}$ we provide the concluding remarks.

 \section{System in two-dimensional self-similar flow}\label{2}
We denote $U=(u-\xi)$ and $V=(v-\eta)$ as pseudo-flow velocity. Then in the self-similar plane $\left(\xi=\dfrac{x}{t}, \eta=\dfrac{y}{t}\right)$, the reduced Euler equations can be written as
\begin{equation}\label{eq: 2.1}
\begin{aligned}
     (\rho U)_\xi+ (\rho V)_\eta+2\rho=0,\\
    U U_\xi+ V U_\eta+ \tau p_\xi+ U=0,\\
        U V_\xi+ V U_\eta+ \tau p_\eta+ V=0.
\end{aligned}
\end{equation}

Here we assume that the flow is irrotational which implies that $v_\xi=u_\eta$. Now we introduce a potential function $\phi(\xi, \eta)$ such that $\phi_\xi=U$ and $\phi_\eta=V$. Then using the last two equations of system \eqref{eq: 2.1} we can easily obtain the pseudo-Bernoulli's law 
\begin{equation}\label{eq: 2.2}
\dfrac{U^2+V^2}{2}+\dfrac{K}{(\tau-b)^\gamma}\left(\dfrac{\gamma+1}{\gamma}+\dfrac{b}{\tau-b}\right)-\dfrac{2a}{\tau}+\phi= a_1
\end{equation}
where $a_1$ is a constant which can be taken as $0$ throughout this article without loss of generality.

The above system $\eqref{eq: 2.1}$ can be reduced into matrix form as
\begin{equation}\label{eq: 2.3}
\begin{aligned}
    \begin{bmatrix}
    u \\ v
    \end{bmatrix}_\xi+\begin{bmatrix}
     \dfrac{-2UV}{c^2-U^2} & \dfrac{c^2-V^2}{c^2-U^2}\\ -1 & 0
    \end{bmatrix} \begin{bmatrix}
    u \\ v
    \end{bmatrix}_\eta=0,
\end{aligned}
\end{equation}
which gives the eigenvalues $\lambda_\pm=\dfrac{UV\pm c\sqrt{U^2+V^2-c^2}}{U^2-c^2}$ with corresponding left eigenvectors $l_\pm=(1, \lambda_\mp)$. Multiplying $l_\pm$ with the system $\eqref{eq: 2.3}$ we obtain the system of characteristic equations as 
\begin{align}\label{eq: 2.4}
\bar{\partial}_{\pm}u+\lambda_\mp \bar{\partial}_{\pm}v=0.
\end{align}

 \subsection{Characteristic equations in terms of characteristic angles}
 Here we provide first order characteristic decompositions of characteristic angles without proof. The proofs of these decompositions can be found in \cite{lai2015expansion}. 
 
 Using \cite{lai2015expansion} we can obtain $U=c\dfrac{\cos\sigma}{\sin \omega}$ and $V=c\dfrac{\sin\sigma}{\sin \omega}$ with
 \begin{align}\label{eq: 2.5}
     \begin{cases}
     \bar{\partial}_-c=\dfrac{\mu^2(\tau)}{\tan \omega}(c\bar{\partial}_-\alpha-2 \sin^2{\omega}),\\
     c \bar{\partial}_-\beta=\Omega(\tau, \omega)\cos^2\omega(c\bar{\partial}_-\alpha-2 \sin^2{\omega})=\dfrac{\Omega(\tau, \omega)}{2 \mu^2(\tau)}\sin{2 \omega}\bar{\partial}_-c, \\
     \bar{\partial}_+c=-\dfrac{\mu^2(\tau)}{\tan \omega}(c\bar{\partial}_+\beta+2 \sin^2{\omega}),\\
     c \bar{\partial}_+\alpha=\Omega(\tau, \omega)\cos^2\omega(c\bar{\partial}_+\beta+2 \sin^2{\omega})=-\dfrac{\Omega(\tau, \omega)}{2 \mu^2(\tau)}\sin{2 \omega}\bar{\partial}_+c, \\
     c \bar{\partial}_\pm \omega=\tan \omega(1+\kappa(\tau) \sin^2{\omega})\bar{\partial}_\pm c+\sin^2{\omega}.
     \end{cases}
 \end{align}
 \subsection{Characteristic decompositions}
 In this section we derive characteristic decomposition form for the variables $\alpha$, $\beta$ and $c$ which is important and very useful for establishing a priori gradient estimates of solution. First we cite the following second order decompositions of $c$ from Lai \cite{lai2015expansion}.
 \begin{p1}\label{p-2.1}
 \textit{The variable $c$ satisfies the following characteristic decompositions}
 \begin{align}\label{eq: 2.6}
 \begin{cases}
     c\bar{\partial}_+\bar{\partial}_- c&=\bar{\partial}_- c\Bigg\{\sin {2\omega}+\dfrac{\bar{\partial}_- c}{2\mu^2(\tau)\cos^2 \omega}+\left(1+\dfrac{\Omega(\tau, \omega) \cos{2 \omega}}{2\mu^2(\tau)}+\tau \kappa'(\tau)\right)\bar{\partial}_+c\Bigg\},\vspace{0.1 cm}\\ 
      c\bar{\partial}_-\bar{\partial}_+ c&=\bar{\partial}_+ c\Bigg\{\sin {2\omega}+\dfrac{\bar{\partial}_+c}{2\mu^2(\tau)\cos^2 \omega}+\left(1+\dfrac{\Omega(\tau, \omega) \cos{2 \omega}}{2\mu^2(\tau)}+\tau \kappa'(\tau)\right)\bar{\partial}_-c\Bigg\}.
      \end{cases}
 \end{align}
 \end{p1}
 Using the above decompositions we are able to prove the following important decompositions for the variable $c$.
 \begin{c1} $c$ satisfies the following second order decompositions
\begin{align}
    \bar{\partial}_+\left(\dfrac{\bar{\partial}_-c}{c}\right)&=\dfrac{\bar{\partial}_-c}{c}\Bigg\{\left(\dfrac{\bar{\partial}_+c+\bar{\partial}_-c}{2 c \mu^2(\tau) \cos^2 \omega}\right)+\dfrac{\sin {2 \omega}}{c}+\left(\tau \kappa'(\tau)-(1+2\kappa(\tau)\sin^2 \omega)\right)\dfrac{\bar{\partial}_+c}{c}\Bigg\}\vspace{0.1 cm},\label{eq: 2.7}\\
 ~~\bar{\partial}_-\left(\dfrac{\bar{\partial}_+c}{c}\right)&=\dfrac{\bar{\partial}_+c}{c}\Bigg\{\left(\dfrac{\bar{\partial}_+c+\bar{\partial}_-c}{2 c \mu^2(\tau) \cos^2 \omega}\right)+\dfrac{\sin {2 \omega}}{c}+\left(\tau \kappa'(\tau)-(1+2\kappa(\tau)\sin^2 \omega)\right)\dfrac{\bar{\partial}_-c}{c}\Bigg\}.  \label{eq: 2.8}
 \end{align}
 \end{c1}
 \begin{proof}
 The proof of this corollary can be obtained by using direct calculations from the decompositions in proposition \ref{p-2.1}. So we omit the details.
 \end{proof}
 \begin{c1}
 $c$ satisfies the following second order equations in homogeneous form
\begin{align} 
    c\bar{\partial}_+ \left(\dfrac{-\bar{\partial}_- c}{\sin^2{\omega}}\right)=\left(\dfrac{-\bar{\partial}_- c}{\sin^2{\omega}}\right)\Bigg\{\left(\dfrac{\bar{\partial}_+ c+\bar{\partial}_- c}{2 \mu^2(\tau) \cos^2 \omega}\right)-\left(\tau \kappa'(\tau)+4 \kappa(\tau) \sin^2{\omega}+2\right)\bar{\partial}_+c\Bigg\},\label{eq: 2.9}\\
      c\bar{\partial}_- \left(\dfrac{\bar{\partial}_+ c}{\sin^2{\omega}}\right)=\left(\dfrac{\bar{\partial}_+ c}{\sin^2{\omega}}\right)\Bigg\{\left(\dfrac{\bar{\partial}_+ c+\bar{\partial}_- c}{2 \mu^2(\tau) \cos^2 \omega}\right)-\left(\tau \kappa'(\tau)+4 \kappa(\tau) \sin^2{\omega}+2\right)\bar{\partial}_-c\Bigg\}. \label{eq: 2.10}
   \end{align}
   \end{c1}
 \begin{proof} This corollary is a direct consequence of the proposition \ref{p-2.1}. Hence we omit its proof.
 \end{proof}
 \begin{p1}\textit{The variables $\alpha$ and $\beta$ satisfy the following second order decompositions}
 \begin{align}
     c\bar{\partial}_+\bar{\partial}_- \alpha+\Psi_1 \bar{\partial}_- \alpha &=\Bigg\{\tan \omega(1-4\sin^2 \omega)+\dfrac{2 \tan \omega \tau \kappa'(\tau)}{\Omega(\tau, \omega)}(\mu^2(\tau))^2\Bigg\}\bar{\partial}_+ \alpha, \label{eq: 2.11} \\
       c\bar{\partial}_-\bar{\partial}_+ \beta+\Psi_2 \bar{\partial}_+ \beta &=\Bigg\{\tan \omega(1-4\sin^2 \omega)+\dfrac{2 \tan \omega \tau \kappa'(\tau)}{\Omega(\tau, \omega)}(\mu^2(\tau))^2\Bigg\}\bar{\partial}_- \beta, \label{eq: 2.12}   
 \end{align}
 in which
 \begin{align*}
     \Psi_1= \sin^2 \omega(2\tan \omega&- \Omega(\tau, \omega) \sin{2 \omega})+\tau \kappa'(\tau)\sin{2 \omega}(\mu^2(\tau))^2-\dfrac{c}{\sin{2 \omega}}\bar{\partial}_-\alpha\\
     &+\Bigg\{\dfrac{1}{\sin{2 \omega}}-\dfrac{\Omega(\tau, \omega) \sin{2\omega}}{2}+\dfrac{\tau \kappa'(\tau)}{\tan \omega}(\mu^2(\tau))^2\Bigg\}c\bar{\partial}_+ \beta,
 \end{align*}
 and
     \begin{align*}
     \Psi_2= \sin^2 \omega(2\tan \omega&- \Omega(\tau, \omega) \sin{2 \omega})+\tau \kappa'(\tau)\sin{2 \omega}(\mu^2(\tau))^2+\dfrac{c}{\sin{2 \omega}}\bar{\partial}_+ \beta\\
     &-\Bigg\{\dfrac{1}{\sin{2 \omega}}-\dfrac{\Omega(\tau, \omega) \sin{2\omega}}{2}+\dfrac{\tau \kappa'(\tau)}{\tan \omega}(\mu^2(\tau))^2\Bigg\}c\bar{\partial}_-\alpha.
     \end{align*}
 \end{p1}
 \begin{proof} Using the decomposition of the variable $c$ from $\eqref{eq: 2.5}$ in $\eqref{eq: 2.6}$ we obtain
 \begin{equation}\label{eq: 2.13}
\begin{aligned}
     c\bar{\partial}_+\Bigg[\dfrac{\mu^2(\tau)c \bar{\partial}_-\alpha}{\tan \omega}-\mu^2(\tau) \sin {2 \omega}\Bigg]=\Bigg[\dfrac{\mu^2(\tau)c \bar{\partial}_-\alpha}{\tan \omega}-\mu^2(\tau) \sin {2 \omega}\Bigg]~~~~~~~~~~~~~~~~~~~~~~~~~~~~~~~\\
     \times \Bigg\{\sin {2 \omega}+\dfrac{c \bar{\partial}_-\beta}{\Omega(\tau, \omega) \cos^2 \omega \sin{2 \omega}}-\left(\dfrac{\Omega(\tau, \omega) \cos{2 \omega}}{2 \mu^2(\tau)}+1+\tau \kappa'(\tau)\right) \dfrac{2 c \mu^2(\tau) \bar{\partial}_+ \alpha}{\Omega(\tau, \omega)\sin{2 \omega}} \Bigg\}
 \end{aligned}
 \end{equation}
 Now we compute
 \begin{equation}\label{eq: 2.14}
 \begin{aligned}
 \bar{\partial}_+\left(\dfrac{\mu^2(\tau)c}{\tan \omega}\right)&=\dfrac{c}{\tan \omega}\left(\mu^2(\tau)\right)' \bar{\partial}_+\tau-\dfrac{\mu^2(\tau)c}{2 \sin^2 \omega}(\bar{\partial}_+\alpha-\bar{\partial}_+\beta)+\dfrac{\mu^2(\tau)}{\tan \omega}\bar{\partial}_+c\\
     &=\dfrac{-\tau (\mu^2(\tau))^3(\kappa^2(\tau))' c \bar{\partial}_+\alpha}{\Omega(\tau, \omega) \tan \omega \sin{2 \omega}}- \mu^2(\tau)-\dfrac{2 (\mu^2(\tau))^2 c \bar{\partial}_+ \alpha}{\Omega(\tau, \omega)\tan \omega \sin{2 \omega}}\\
     &~~~~+\dfrac{\mu^2(\tau)c \bar{\partial}_+ \alpha}{\Omega(\tau, \omega) \sin{2 \omega} \sin^2{\omega}}\left(\tan \omega- \dfrac{\Omega(\tau, \omega) \sin {2 \omega}}{2}\right)
     \end{aligned}
 \end{equation}
 and \begin{equation}\label{eq: 2.15}
     \begin{aligned}
     c\bar{\partial}_+\left(\mu^2(\tau)\sin {2\omega} \right)&= c \sin{2 \omega} (\mu^2(\tau))' \bar{\partial}_+\tau+ \mu^2(\tau) c \cos{2 \omega}[\bar{\partial}_+\alpha -\bar{\partial}_+\beta]\\
     &=\dfrac{-\tau (\mu^2(\tau))^3(\kappa^2(\tau))' c \bar{\partial}_+\alpha}{\Omega(\tau, \omega)}- 2 \mu^2(\tau) \sin^2\omega+\mu^2(\tau) \sin^2(2 \omega )\\
    &~~~-\dfrac{2 \mu^2(\tau) \cos{2 \omega} c \bar{\partial}_+\alpha}{\Omega(\tau, \omega)\sin{2 \omega}}\left(\tan \omega - \dfrac{\Omega(\tau, \omega)\sin {2 \omega}}{2}\right).
     \end{aligned}
 \end{equation}
  so that the L.H.S. of $\eqref{eq: 2.13}$ becomes
 \begin{equation}\label{eq: 2.16}
     \begin{aligned}
     c\bar{\partial}_+\Bigg[\dfrac{\mu^2(\tau)c\bar{\partial}_-\alpha}{\tan \omega}&-\mu^2(\tau)\sin {2\omega}\Bigg]= \dfrac{\mu^2(\tau)c}{\tan \omega}c \bar{\partial}_+\bar{\partial}_-\alpha+ c \bar{\partial}_-\alpha \bar{\partial}_+\left(\dfrac{\mu^2(\tau)c}{\tan \omega}\right)-c \bar{\partial}_+\left(\mu^2(\tau)\sin {2\omega}\right)\\
     &=\dfrac{\mu^2(\tau)c}{\tan \omega}c \bar{\partial}_+\bar{\partial}_-\alpha
     +\dfrac{\tau (\mu^2(\tau))^3(\kappa^2(\tau))' c \bar{\partial}_+\alpha}{\Omega(\tau, \omega)}+ 2 \mu^2(\tau) \sin^2\omega-\mu^2(\tau) \sin^2 (2 \omega)\\
     &+ c \bar{\partial}_-\alpha\Bigg\{\dfrac{\mu^2(\tau)c \bar{\partial}_+ \alpha}{\Omega(\tau, \omega) \sin{2 \omega} \sin^2{\omega}}\left(\tan \omega
     - \dfrac{\Omega(\tau, \omega) \sin {2 \omega}}{2}\right)-\dfrac{\tau (\mu^2(\tau))^3 (\kappa^2(\tau))' c \bar{\partial}_+\alpha}{\Omega(\tau, \omega) \tan \omega \sin{2 \omega}}\\
     &-\dfrac{2 (\mu^2(\tau))^2 c \bar{\partial}_+ \alpha}{\Omega(\tau, \omega)\tan \omega \sin{2 \omega}}- \mu^2(\tau)\Bigg\}+\dfrac{2 \mu^2(\tau) \cos{2 \omega} c \bar{\partial}_+\alpha}{\Omega(\tau, \omega)\sin{2 \omega}}\left(\tan \omega - \dfrac{\Omega(\tau, \omega)\sin {2 \omega}}{2}\right).
     \end{aligned}
 \end{equation}
While the R.H.S. of $\eqref{eq: 2.13}$ is 
 \begin{equation}\label{eq: 2.17}
     \begin{aligned}
\dfrac{\mu^2(\tau)c}{\tan \omega} \bar{\partial}_-\alpha 
\Bigg\{
    \sin{2 \omega}-2 \tan \omega+\dfrac{c\bar{\partial}_-\alpha}{\sin{2 \omega}}-\left(1+\dfrac{\Omega(\tau, \omega) \cos{2 \omega}}{2 \mu^2(\tau)}+\tau \kappa'(\tau)\right) \dfrac{2 c \mu^2(\tau) \bar{\partial}_+ \alpha}{\Omega(\tau, \omega)\sin{2 \omega}}\Bigg\}\\+\left(1+\dfrac{\Omega(\tau, \omega) \cos{2 \omega}}{2 \mu^2(\tau)}+\tau \kappa'(\tau)\right) \dfrac{2 (\mu^2(\tau))^2 c \bar{\partial}_+ \alpha}{\Omega(\tau, \omega)}
+ 2 \mu^2(\tau) \sin^2\omega-\mu^2(\tau) \sin^2(2 \omega).
     \end{aligned}
 \end{equation}
 Comparing L.H.S. and R.H.S. of $\eqref{eq: 2.13}$ we have
 \begin{equation}\label{eq: 2.18}
     \begin{aligned}
c \bar{\partial}_+\bar{\partial}_-\alpha+ c \bar{\partial}_-\alpha\Bigg\{\dfrac{\tan \omega \bar{\partial}_+ \alpha}{\Omega(\tau, \omega) \sin{2 \omega} \sin^2{\omega}}\left(\tan \omega
     - \dfrac{\Omega(\tau, \omega) \sin {2 \omega}}{2}\right)-\dfrac{2 \mu^2(\tau) \bar{\partial}_+ \alpha}{\Omega(\tau, \omega) \sin{2 \omega}}\\
    - \dfrac{\tan \omega}{c}-\dfrac{\tau ((\mu^2(\tau))^2 (\kappa^2(\tau))'  \bar{\partial}_+\alpha}{\Omega(\tau, \omega)\sin{2 \omega}}\Bigg\}\\
     +\dfrac{\tau (\mu^2(\tau))^2 (\kappa^2(\tau))' \tan \omega \bar{\partial}_+\alpha}{ \Omega(\tau, \omega)}+\dfrac{2 \tan \omega \cos{2 \omega} \bar{\partial}_+\alpha}{\Omega(\tau, \omega)\sin{2 \omega}}\left(\tan \omega - \dfrac{\Omega(\tau, \omega)\sin {2 \omega}}{2}\right)\\
=\bar{\partial}_-\alpha \Bigg\{
    \sin{2 \omega}-2 \tan \omega +\dfrac{c\bar{\partial}_-\alpha}{\sin{2 \omega}}-\left(1+\dfrac{\Omega(\tau, \omega) \cos{2 \omega}}{2 \mu^2(\tau)}+\tau \kappa'(\tau)\right) \dfrac{2 c \mu^2(\tau) \bar{\partial}_+ \alpha}{\Omega(\tau, \omega)\sin{2 \omega}}\Bigg\}\\
    +\left(1+\dfrac{\Omega(\tau, \omega) \cos{2 \omega}}{2 \mu^2(\tau)}+\tau \kappa'(\tau)\right) \dfrac{2 \mu^2(\tau) \tan \omega \bar{\partial}_+ \alpha}{\Omega(\tau, \omega)},
     \end{aligned}
 \end{equation}
 which can be reduced as 
 \begin{equation}\label{eq: 2.19}
     \begin{aligned}
     c\bar{\partial}_+\bar{\partial}_- \alpha+\Psi_1 \bar{\partial}_- \alpha =\Upsilon_1 \bar{\partial}_+ \alpha\vspace{0.2 cm}\\
     \end{aligned}
 \end{equation}
 where the coefficient 
 \begin{equation}\label{eq: 2.20}
 \begin{aligned}
 \Psi_1&=\left(\dfrac{\sin {2 \omega}}{2}+\dfrac{c \bar{\partial}_+ \beta}{2 \tan \omega}\right)\Bigg[\dfrac{\tan \omega}{\sin^2 \omega}\left(\tan \omega-\dfrac{\Omega(\tau, \omega) \sin{2 \omega}}{2}\right)+\Omega(\tau, \omega)\cos{2 \omega}+2 \tau \kappa'(\tau)(\mu^2(\tau))^2\Bigg]\\
     &~~~~~~~~~~~~~~~~~~~~~~~~~~~~~~~~~~~~~~~~~~~~~~~~~~~~~~~~~~~~~~~~~~~~~~~~~~~~~~~~~~~~+\tan \omega-\sin {2 \omega}-\dfrac{c\bar{\partial}_-\alpha}{\sin{2 \omega}}.
 \end{aligned}
 \end{equation}
 By a direct calculation, we obtain the coefficient of $c \bar{\partial}_+\beta$ as\\
 \begin{equation}\label{eq: 2.21}
     \begin{aligned}
\dfrac{1}{\sin {2 \omega}}-\dfrac{\Omega(\tau, \omega) \sin{2 \omega}}{2}+\dfrac{\tau \kappa'(\tau)}{\tan \omega}(\mu^2(\tau))^2.
     \end{aligned}
 \end{equation}
Also, we calculate directly\\
\begin{equation}\label{eq: 2.22}
    \begin{aligned}
\dfrac{\sin {2 \omega}}{2}\Bigg\{\dfrac{\tan \omega}{\sin^2 \omega}\left(\tan \omega-\dfrac{\Omega(\tau, \omega) \sin{2 \omega}}{2}\right)+\Omega(\tau, \omega)\cos{2 \omega}+2  \tau \kappa'(\tau)(\mu^2(\tau))^2\Bigg\}\\
+\tan \omega-\sin{2 \omega}- \dfrac{c\bar{\partial}_-\alpha}{\sin{2 \omega}}\\  
=\sin^2{\omega}(2 \tan \omega-\Omega(\tau, \omega)\sin{2 \omega})+\sin{2 \omega} \tau \kappa'(\tau)(\mu^2(\tau))^2- \dfrac{c\bar{\partial}_-\alpha}{\sin{2 \omega}}, 
    \end{aligned}
\end{equation}
Using $\eqref{eq: 2.21}$ and $\eqref{eq: 2.22}$ in $\eqref{eq: 2.20}$ we obtain the value of $\Psi_1$.

Furthermore, the coefficient
\begin{equation}\label{eq: 2.23}
 \begin{aligned}
 \Upsilon_1= \dfrac{-2 \tan \omega \cos{2 \omega} }{\Omega(\tau, \omega)\sin{2 \omega}}\left(\tan \omega - \dfrac{\Omega(\tau, \omega)\sin {2 \omega}}{2}\right)+\left(1+\dfrac{\Omega(\tau, \omega) \cos{2 \omega}}{2 \mu^2(\tau)}\right) \dfrac{2 \mu^2(\tau) \tan \omega}{\Omega(\tau, \omega)}\\
 +\dfrac{2 \tan \omega \tau \kappa'(\tau)}{\Omega(\tau, \omega)}(\mu^2(\tau))^2
 \end{aligned}
\end{equation}
By straightforward computation we have
\begin{equation}\label{eq: 2.24}
    \begin{aligned}
&\dfrac{-2 \tan \omega \cos{2 \omega} }{\Omega(\tau, \omega)\sin{2 \omega}}\left(\tan \omega - \dfrac{\Omega(\tau, \omega)\sin {2 \omega}}{2}\right)+\left(1+\dfrac{\Omega(\tau, \omega) \cos{2 \omega}}{2 \mu^2(\tau)}\right) \dfrac{2 \mu^2(\tau) \tan \omega}{\Omega(\tau, \omega)}\\    
&=\dfrac{-2 \tan \omega}{\Omega(\tau, \omega)\sin{2 \omega}}\Big[\mu^2(\tau)\sin{2 \omega}-\Omega(\tau, \omega)\sin{2 \omega}\cos{2 \omega}-\cos{2 \omega} \tan{\omega}\Big]\\
&=2\tan \omega\left(\dfrac{1}{2}-\cos{2 \omega}\right)\\
&=\tan \omega(1-4 \sin^2{\omega}).
    \end{aligned}
\end{equation}
Then the use of $\eqref{eq: 2.24}$ in $\eqref{eq: 2.23}$ yields the value of $\Upsilon_1$.

Using the values of $\Psi_1$ and $\Upsilon_1$ in $\eqref{eq: 2.19}$ we can obtain the proof of the first part of the Lemma. In a similar manner the other part of the Lemma can be proved.
\end{proof}
 \section{Formulation of the problem and local solution}\label{3}
In this section we give a brief description of our problem for $\eqref{eq: 1}$. Let us assume that $R_{14}(\eta)$ be the planar rarefaction wave for $\eqref{eq: 1}$ connecting two constant states $(\rho_1, 0, v_1)$ and $(\rho_4, 0, 0)$ ($c_1>c_4>0$, $v_1>0$) in the self-similar plane which is defined as
\begin{align}{\label{eq: 3.1}}
R_{14}(\eta)=
\begin{cases}
    \eta=c+v(\rho),~~~~~~~~~~~~ \eta_4\leq\eta\leq \eta_1\\
    v=\displaystyle \int_{\rho_4}^{\rho}\dfrac{c}{\rho}d\rho,~~~~~~~~~~~~~ 0<\rho_4<\rho<\rho_1\\
    u=u_1=u_4=0,
\end{cases}
\end{align}
in the region $\xi>0$, where $\eta_i=v_i+c_i$, $i=1, 4$.

The value $v_1$ is obtained from the solution and we denote $A$ by $(0, v_1+c_1)$ and $B$ as the intersection point of positive characteristic passing through $A$ and the bottom boundary in the rarefaction wave region. Now under these assumptions we define our problem as follows.
\subsection{Problem} 
\begin{figure}
    \centering
    \includegraphics[width=5.5 in]{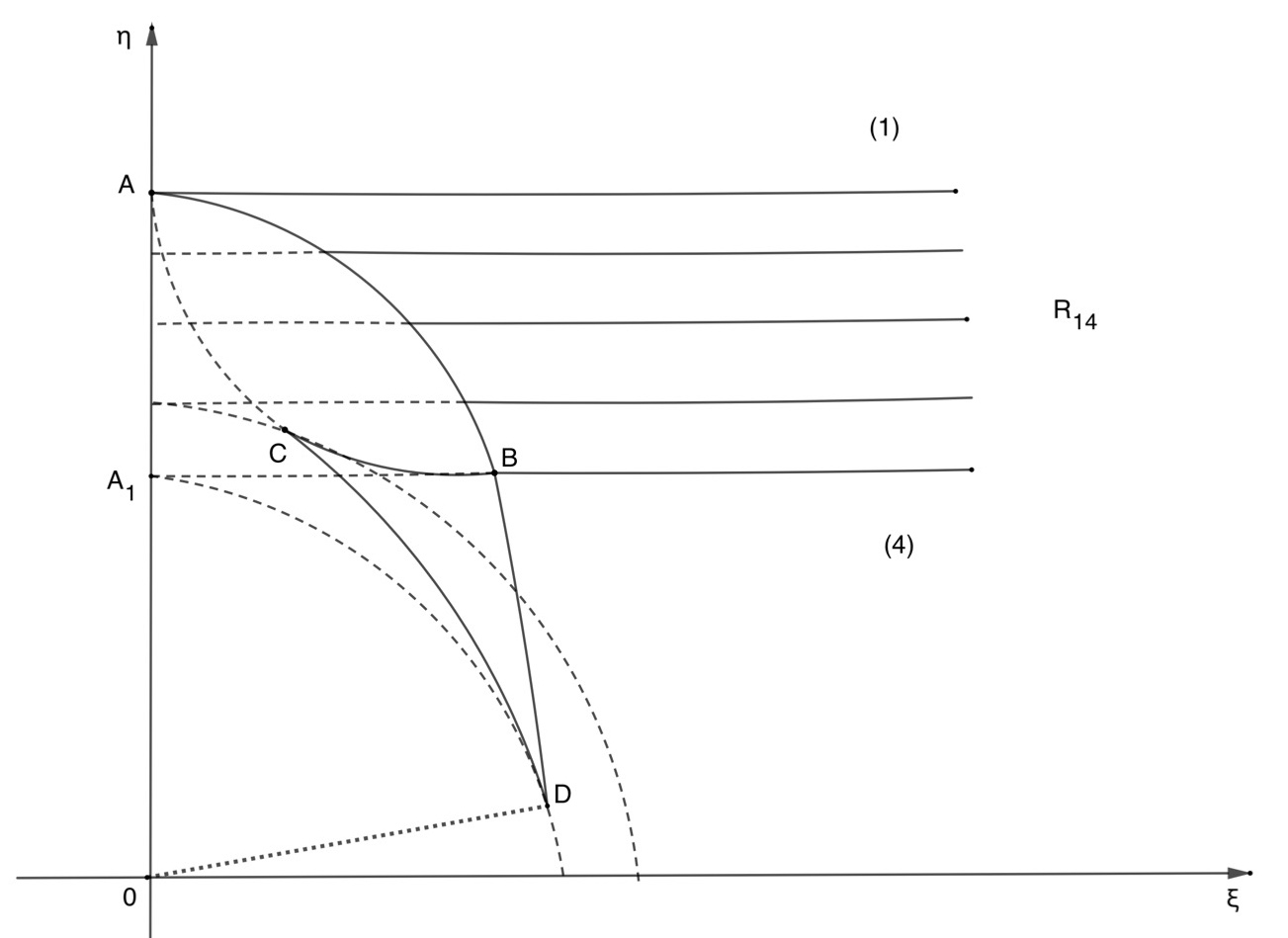}
    \caption{The semi-hyperbolic patch.}
    \label{fig:my_label1}
\end{figure}

\textit{Let us consider a positive characteristic $\overline{AB}$ in a planar rarefaction wave region. Let $\overline{BD}$ be the tangential extension of $\overline{AB}$ into a constant state such that both the points $A$ and $D$ are sonic. Let $\overline{BC}$ be a strictly convex negative characteristic such that the endpoint $C$ is a sonic point. Then establish a solution in a maximal hyperbolic region with the boundary data combined on the curve $\overline{ABD}$ and $\overline{BC}$ which means that the region starts from sonic points and ends on either a sonic curve or an envelope of the positive family of characteristics. Further, establish the regularity of solution; see Figure $\ref{fig:my_label1}$.} 

\subsection{Estimates of boundary data}
The boundary data on the characteristics $\overline{AB}$ and $\overline{BC}$ are \cite{li2011semi}
\begin{equation}{\label{eq: 3.2}}
\begin{cases}
~~~~~~~~~~~~\beta|_{\overline{AB}}=0,~~~~~~~\pi/2\leq \alpha|_{\overline{BC}}\leq \pi+\beta_C,
\vspace{0.2 cm}\\
\pi/2\leq \alpha|_{\overline{AB}}\leq \pi, ~~~~~~~\beta_C\leq \beta|_{\overline{BC}}\leq 0,\vspace{0.2 cm}\\
~-\pi/2<\beta_C<0, ~~~~~~~\dfrac{\pi}{4}<\dfrac{\alpha-\beta}{2}<\dfrac{\pi}{2},\vspace{0.2 cm}\\
~~~~~~~~\bar{\partial}_+ \alpha|_{\overline{BA}}>0,~~~~~~
        ~~\bar{\partial}_+ c|_{\overline{BA}}>0,\vspace{0.2 cm}\\
        ~~~~~-\bar{\partial}_- \beta|_{\overline{BC}}<0,  ~~~~~-\bar{\partial}_- c|_{\overline{BC}}>0.
\end{cases}
\end{equation}
where $\beta_C$ is the inclination angle of the negative characteristic $\overline{BC}$ at point $C$. 

Hence our problem is a Goursat problem where the boundaries $\overline{AB}$ and $\overline{BC}$ are characteristic boundaries starting from the point $B$.

\subsection{Existence of local solution}
To prove the existence of local solution, we parameterize the variables $\alpha$, $\beta$ and $c$ on the positive and negative characteristics $\overline{AB}$ and $\overline{BC}$ as a function of parameter $s$. Further, we assume the point $B$ as the origin ($s=0$). Then using $\eqref{eq: 2.5}$ and $\eqref{eq: 3.2}$ we obtain the boundary data in the form of parameter $s$ as follows:

On $\overline{AB}$ we have
\begin{equation}{\label{eq: 3.3}}
    \begin{aligned}
    \begin{cases}
\dfrac{\pi}{2}< \alpha(s) < \pi,~~~~\beta(s)=0, \vspace{0.15 cm}   \\ 
\dfrac{\pi}{4}<\dfrac{\alpha(s)-\beta(s)}{2}< \dfrac{\pi}{2}, \vspace{0.15 cm} \\
\dfrac{d \alpha}{ds}>0,~~~~\dfrac{d c}{ds}>0, \vspace{0.15 cm}\\
c\dfrac{d \alpha}{ds}=-\dfrac{\sin{2 \omega}}{2 \mu^2(\tau)}\Omega(\tau, \omega) \dfrac{dc}{ds}, \vspace{0.15 cm}\\
c\dfrac{d \beta}{ds}=-\dfrac{\tan{\omega}}{\mu^2(\tau)} \left(\dfrac{dc}{ds}\right)-2 \left(\sqrt{\left(\dfrac{d \xi}{ds}\right)^2+\left(\dfrac{d \eta}{ds}\right)^2}\right)\sin^2{\omega}
\end{cases}
\end{aligned}
\end{equation}
and on $\overline{BC}$ we have
\begin{equation}{\label{eq: 3.4}}
    \begin{aligned}
\begin{cases}
\dfrac{\pi}{2}< \alpha(s) <\pi+\beta_C,~~~~\beta_C\leq\beta(s) \leq 0, \vspace{0.15 cm}\\
\dfrac{\pi}{4}<\dfrac{\alpha(s)-\beta(s)}{2}<\dfrac{\pi}{2}, \vspace{0.15 cm}\\
\dfrac{d \beta}{ds}<0, ~~~~\dfrac{d c}{ds}>0, \vspace{0.15 cm}\\
c\dfrac{d \beta}{ds}=\dfrac{\sin{2 \omega}}{2 \mu^2(\tau)}\Omega(\tau, \omega) \dfrac{dc}{ds}, \vspace{0.15 cm}\\
c\dfrac{d \alpha}{ds}=\dfrac{\tan{\omega}}{\mu^2(\tau)} \left(\dfrac{dc}{ds}\right)+2 \left(\sqrt{\left(\dfrac{d \xi}{ds}\right)^2+\left(\dfrac{d \eta}{ds}\right)^2}\right)\sin^2{\omega}.
\end{cases}
    \end{aligned}
\end{equation}

Using these boundary values we can prove the existence of local solution for our Goursat problem. We summarize this in the following Lemma.
\begin{l1}\label{l-3.1}
\textit{(Local solution) For sufficiently small $\epsilon>0$ there exists a unique $C^1$ solution for the Goursat problem $\eqref{eq: 2.1}$, $\eqref{eq: 3.2}$ in a small domain $D_\epsilon$ closed by the boundaries $\overline{AB}$ and $\overline{BC}$ and a level curve $\omega=\epsilon$. Further, this solution satisfies
\begin{align}\label{eq: 3.5}
\dfrac{\pi}{4}< \omega< \dfrac{\pi}{2},~~ \pm \bar{\partial}_\pm c>0,~~~ \pm \bar{\partial}_\pm \alpha>0,~ ~~\pm \bar{\partial}_\pm\beta<0~~~ in~~D_\epsilon.  
\end{align}}
\end{l1}
\begin{proof}
It is observed from $\eqref{eq: 3.3}$ and $\eqref{eq: 3.4}$ that the compatibility conditions hold at the point B. Therefore, using the method of characteristics \cite{li1985boundary}, we conclude that the Goursat problem $\eqref{eq: 2.1}$, $\eqref{eq: 3.2}$ admits a unique $C^1$ local solution. From boundary estimate $\eqref{eq: 3.2}$ we have 
\begin{figure}
    \centering
    \includegraphics[width= 3.5 in]{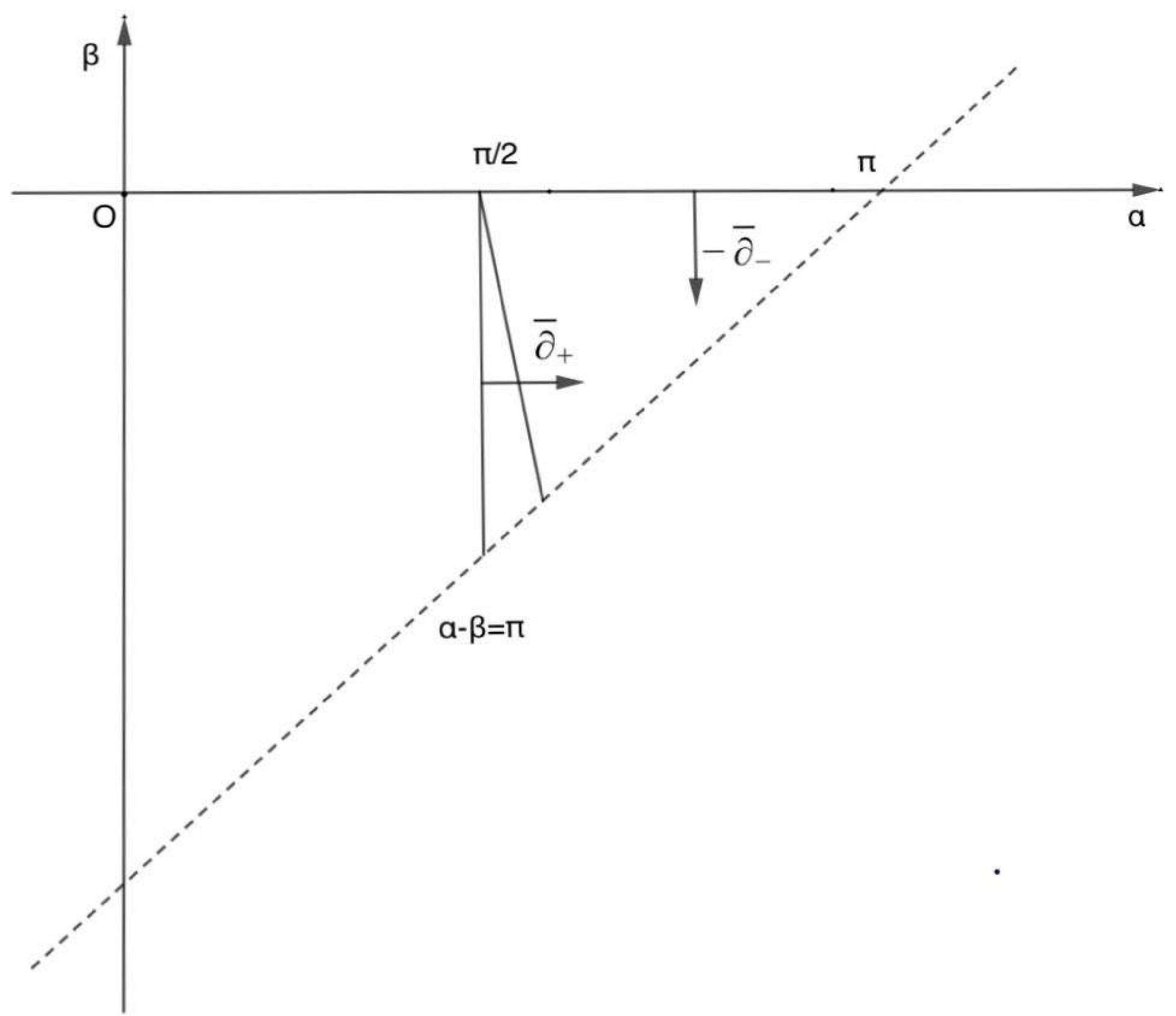}
    \caption{Invariant triangle}
    \label{fig:my_label2}
\end{figure}

\begin{center}
$\bar{\partial}_+ c|_{\overline{AB}}>0$ and $-\bar{\partial}_- c|_{\overline{BC}}>0$.
\end{center}
So using the characteristic decompositions from $\eqref{eq: 2.6}$ we can conclude that
\begin{center}
$\bar{\partial}_+c>0$ and $-\bar{\partial}_-c>0$ in $D_\epsilon$.
\end{center}
Now we prove $\omega> \frac{\pi}{4}$ in $D_\epsilon$ using method of contradiction. Let us assume that there exists a point $G$ such that $\omega_{G}=\frac{\pi}{4}$. Then using $\eqref{eq: 2.5}$ we observe that $\bar{\partial}_+ \omega|_{G}>0$ which provides a contradiction, since $\omega|_{\overline{AB}}>\frac{\pi}{4}$. Hence we have $\omega>\frac{\pi}{4}$ in $D_{\epsilon}$.

Further, the fact $\omega>\frac{\pi}{4}$ and $m(\tau)<1$ shows that $\Omega(\tau, \omega)<0$ in $D_\epsilon$. Hence, by characteristic decompositions $\eqref{eq: 2.5}$ we have 
\begin{center}
$\bar{\partial}_+\alpha>0,~~~\bar{\partial}_+\beta<0,~~ -\bar{\partial}_-\beta<0$ in $D_{\epsilon}$.
\end{center}
From the boundary data $\eqref{eq: 3.2}$ we have $-\bar{\partial}_-\alpha|_{\overline{BC}}>0$. So using characteristic decomposition $\eqref{eq: 2.11}$ we obtain $-\bar{\partial}_-\alpha>0$ in the domain $D_{\epsilon}$.

Therefore, using $\bar{\partial}_+\alpha>0$ and $-\bar{\partial}_-\beta<0$ we obtain the invariant triangle, see Figure $\ref{fig:my_label2}$ (For more details on invariant regions, see \cite{smoller2012shock})
$$\alpha>\pi/2, ~~~\beta<0,~~~\frac{\pi}{4}<\omega<\frac{\pi}{2}.$$ 
Hence the Lemma is proved.
\end{proof}
\begin{l1}\label{l-3.2}
\textit{If the Goursat problem $\eqref{eq: 2.1}$, $\eqref{eq: 3.2}$ admits a unique $C^1$ solution in the domain $D_\epsilon, \epsilon \in (\pi/4, \pi/2)$. Then there exist a positive constant $\mathcal{H}$ such that
\begin{align}
    ||(D\alpha, D\beta, Dc)||_{D_{\epsilon}}<\mathcal{H}.
\end{align}}
\end{l1}
\begin{proof} Using the Lemma $\ref{l-3.1}$, $\eqref{eq: 2.5}$, relations $\sin \beta \bar{\partial}_+-\sin \alpha \bar{\partial}_-=-\sin{2\omega}\partial_\xi$ and $\cos \beta \bar{\partial}_+-\cos \alpha \bar{\partial}_-=\sin{2 \omega}\partial_\eta$, we can easily prove this Lemma.
\end{proof}
\section{Uniform bounds of $\pm \bar{\partial}_\pm c$}\label{4}
In this section, we give uniform lower and upper bounds of $ \bar{\partial}_+ c$ and $-\bar{\partial}_- c$ which are useful for further analysis in the succeeding sections.
\subsection{Upper bounds of $\pm \bar{\partial}_\pm c$}
We provide the upper bounds of $\bar{\partial}_+ c$ and $-\bar{\partial}_- c$ in the following Lemma:
\begin{l1}\label{l-4.1}
\textit{If the Goursat problem $\eqref{eq: 2.1}$, $\eqref{eq: 3.2}$ admits a $C^1$ solution in the domain $D_{\epsilon}$ where $\epsilon \in\left(\frac{\pi}{4}, \frac{\pi}{2}\right)$, then we have
\begin{align}\label{eq: 4.1}
    \left(\dfrac{\bar{\partial}_+c}{\sin^2{\omega}}, -\dfrac{\bar{\partial}_-c}{\sin^2{\omega}}\right)\in (0, M_{\epsilon}]\times (0, M_{\epsilon}]
\end{align}
where $M_{\epsilon}=\max \Bigg\{\underset{\overline{BA_\epsilon}}{\sup} \left(\dfrac{\bar{\partial}_+c}{\sin^2{\omega}}\right),~ \underset{\overline{BC_\epsilon}}{\sup} \left(\dfrac{-\bar{\partial}_-c}{\sin^2{\omega}}\right)$ \Bigg\} such that $A_{\epsilon}$ and $C_{\epsilon}$ are points on $AB$ and $BC$, respectively with $\omega(A_{\epsilon})=\omega(C_{\epsilon})=\epsilon$.}
\end{l1}
\begin{proof}
In order to prove this Lemma we need to prove that for any $x>0$ we have $$\left(\dfrac{\bar{\partial}_+c}{\sin^2{\omega}}, -\dfrac{\bar{\partial}_-c}{\sin^2{\omega}}\right)\in (0, M_{\epsilon}+x)\times (0, M_{\epsilon}+x)~in~ D_{\epsilon}.$$
In the contrary let us assume that the Lemma is not valid. Then we must have at least one point $E$ in $D_{\epsilon}$ such that $\left(\frac{\bar{\partial}_+c(E)}{\sin^2{\omega}}, -\frac{\bar{\partial}_-c(E)}{\sin^2{\omega}}\right) \in \bigcup\limits_{i=1}^2 \Gamma_i$ and  $\left(\frac{\bar{\partial}_+c}{\sin^2{\omega}}, -\frac{\bar{\partial}_-c}{\sin^2{\omega}}\right)\in (0, M_{\epsilon}+x)\times (0, M_{\epsilon}+x)$ for any $(\xi, \eta)\in \Sigma_{E}/ \{E\}$ where $\Gamma_1=(0, M_{\epsilon}+x]\times \{M_{\epsilon}+x\}$ and $\Gamma_2=\{M_{\epsilon}+x\}\times (0, M_{\epsilon}+x]$ with $\Sigma_{E}$ being the region bounded by $\overline{BE_{1}}$, $\overline{BE_{2}}$, $\overline{EE_{1}}$ and $\overline{EE_{2}}$ such that the point $E_{1}$ is the point of intersection of the boundary $\overline{AB}$ and the negative characteristic passing through the point $E$ and similarly the point $E_{2}$ is the point of intersection of the boundary $\overline{BC}$ and the positive characteristic passing through the point $E$; see Figure $\ref{fig: my_label3}$.
\begin{figure}
    \centering
   \includegraphics[width= 6.1 in]{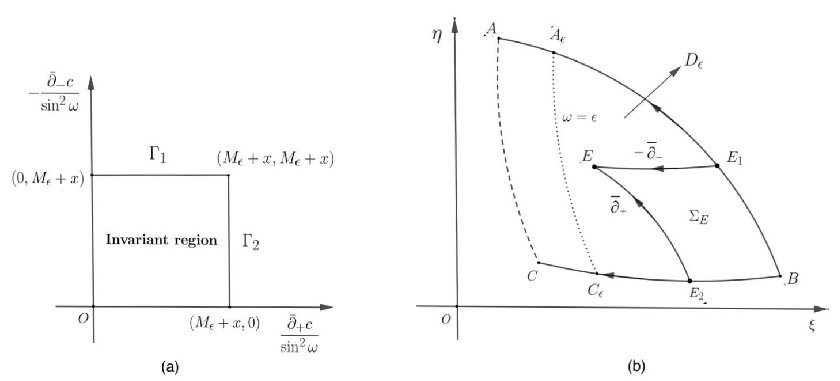}
    \caption{(a) Invariant region of $\left(\dfrac{\bar{\partial}_+c}{\sin^2{\omega}}, -\dfrac{\bar{\partial}_-c}{\sin^2{\omega}}\right)$; (b) Domain $\Sigma_{E}$ }
    \label{fig: my_label3}
\end{figure}
Suppose that $\left(\frac{\bar{\partial}_+c(E)}{\sin^2{\omega}}, -\frac{\bar{\partial}_-c(E)}{\sin^2{\omega}}\right) \in \Gamma_1$, then using $\eqref{eq: 2.9}$ we have $\bar{\partial}_+\left(\frac{-\bar{\partial}_-c(E)}{\sin^2{\omega}}\right)<0$ which leads to a contradiction since $\frac{-\bar{\partial}_-c(E)}{\sin^2{\omega}}=\{M_{\epsilon}+x\}$. Similarly, if $\left(\frac{\bar{\partial}_+c(E)}{\sin^2{\omega}}, -\frac{\bar{\partial}_-c(E)}{\sin^2{\omega}}\right) \in \Gamma_2$, then again by exploiting $\eqref{eq: 2.10}$, we have $\bar{\partial}_-\left(\frac{\bar{\partial}_+c(E)}{\sin^2{\omega}}\right)<0$ which leads to a contradiction since $\frac{\bar{\partial}_+c(E)}{\sin^2{\omega}}=\{M_{\epsilon}+x\}$. Therefore, our assumption is wrong which proves the Lemma.
\end{proof}
\subsection{Lower bounds of $\pm \bar{\partial}_\pm c$}
In the following Lemma we find the lower bounds of $\dfrac{\bar{\partial}_+ c}{c}$ and $-\dfrac{\bar{\partial}_- c}{c}$ which directly provides us lower bounds of $\bar{\partial}_+c$ and $-\bar{\partial}_-c$.
\begin{l1}\label{l-4.2}
\textit{If the Goursat problem $\eqref{eq: 2.1}$, $\eqref{eq: 3.2}$ admits a $C^1$ solution in the domain $D_{\epsilon}$ then the functions $\frac{\bar{\partial}_+c}{c}$ and $\frac{\bar{\partial}_-c}{c}$ satisfy
\begin{center}
   $0 < \overline{m} e^{-\Hat{\kappa} d}\leq \dfrac{\bar{\partial}_+c}{c},~~~   0< \bar{m} e^{-\Hat{\kappa} d}\leq -\dfrac{\bar{\partial}_-c}{c}$
\end{center}
where $2 \overline{m}=\min \{\underset{\overline{BA_*}}{\min} \frac{\bar{\partial}_+c}{c}, \underset{\overline{BC_*}}{\min} (-\frac{\bar{\partial}_-c}{c})\}>0$, $d$ is the diameter of the domain $D_{\epsilon}$, $\Hat{\kappa}(\tau)=\dfrac{3+4 \kappa(\tau)}{2 c_4}$, $A_*$ is the point of intersection of negative characteristic and the boundary $\overline{BA}$ and $C_*$ is the point of intersection of positive characteristic and the boundary $\overline{BC}$ such that both the characteristic curves start from a single point $P$ in the domain $D_{\epsilon}$.}
\begin{figure}
    \centering
    \includegraphics[width=5 in]{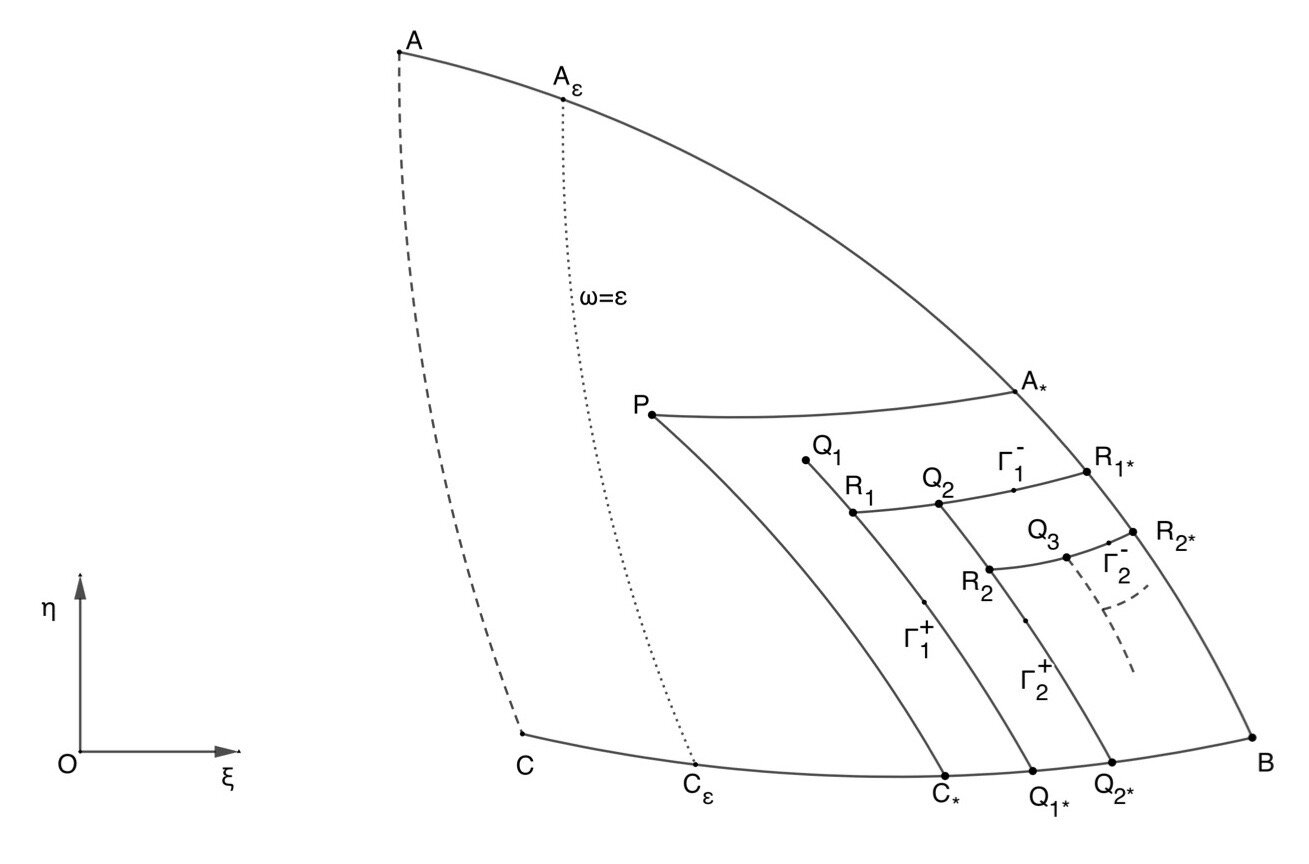}
    \caption{Characteristic curves passing through different points in $D_\epsilon$}
    \label{fig:my_label4}
\end{figure}
\end{l1}
\begin{proof}
The uniform lower bounds of $\dfrac{\bar{\partial}_+c}{c}$ and $-\dfrac{\bar{\partial}_-c}{c}$ are established using the new variables
\begin{align}\label{eq: 4.2}
    \overline{R}=A(\kappa) \sin{\omega}\sqrt{1+\kappa(\tau) \sin^2{\omega}}~\dfrac{\bar{\partial}_+c}{c}>0,~~~~ \overline{S}=-A(\kappa) \sin{\omega}\sqrt{1+\kappa(\tau) \sin^2{\omega}}~\dfrac{\bar{\partial}_-c}{c}>0
\end{align}
where $A(\kappa(\tau))=\kappa(\tau) \exp\left(\displaystyle \int \dfrac{-\kappa'(\tau)\sin^2{\omega}}{2(1+\kappa(\tau) \sin^2{\omega})} d \tau \right)>0$ in the region $D_\epsilon$.\\            

Then using $\eqref{eq: 2.7}$ and $\eqref{eq: 2.8}$ we obtain
\begin{align}\label{eq: 4.3}
\bar{\partial}_- \overline{R}= \overline{R} \varPi,~~~~~~~\bar{\partial}_+ \overline{S}= \overline{S} \varPi
\end{align}
where 
\begin{align}\label{eq: 4.4}
\varPi=\dfrac{(\overline{R}-\overline{S})(1+\kappa(\tau))}{A(\kappa(\tau)) \sin{2\omega} \cos{\omega}\sqrt{1+\kappa(\tau) \sin^2{\omega}}}+\dfrac{\sin{2 \omega}}{2c}\left(\dfrac{3+4\kappa(\tau) \sin^2{\omega}}{1+\kappa(\tau) \sin^2{\omega}}\right).
\end{align}
Now we consider the following two cases. \vspace{0.2 cm}\\
\textbf{Case A:} Let us assume that $\overline{S}\geq \overline{R}$ holds entirely in the region $PA_*BC_*P$ for each point $P$.\vspace{0.2 cm}\\
Then using $\eqref{eq: 4.3}$ we get
\begin{align}\label{eq: 4.5}
\bar{\partial}_- \ln \overline{R} \leq\dfrac{3+4\kappa(\tau)}{2c} \leq\dfrac{3+4\kappa(\tau)}{2c_4}=\Hat{\kappa}(\tau).
\end{align}
Then integrating $\eqref{eq: 4.5}$ along the negative characteristic from $A_*$ to $P$ yields
$$\overline{S}_P>\overline{R}_P\geq \overline{R}_{A_*}e^{-\Hat{\kappa}(\eta_P-\eta_{A_*})}\geq \overline{m}e^{-\Hat{\kappa} d}.$$
\textbf{Case B:} Let us assume that there exists a point $Q_1$ in $PA_*BC_*P$ such that $\overline{S}<\overline{R}$ at $Q_1$. Now we draw a positive characteristic curve $\Gamma_1^+$ starting from $Q_1$ to $Q_{1*}$ which lies on the boundary $\overline{BC_*}$.
 
 If $\varPi\geq 0$ for all the points on $\Gamma_1^+$ then using $\eqref{eq: 4.3}$ we have 
  \begin{align*}
\bar{\partial}_+\overline{S}>0~~~~~~~~~ on~~~~~~ \Gamma_1^+
  \end{align*}
which implies that
 \begin{align*}
\overline{R}_{Q_1}>\overline{S}_{Q_1}\geq \overline{S}_{Q_{1*}}\geq \overline{m}\geq \overline{m}e^{-\Hat{\kappa} d}.
 \end{align*}
Otherwise there exists a point on $\Gamma_1^+$ such that $\varPi<0$ at that point. But $\varPi>0$ at $Q_1$ so we use continuity of $\varPi$ to see that $\varPi> 0$ holds for all points in $N_1\cap \Gamma_1^+$ for some neighbourhood $N_1$ of $Q_1$ and $\varPi=0$ at $R_1=\partial N_1\cap \Gamma_1^+$.

So we see that
\begin{center}
  $\bar{\partial}_+\overline{S}>0$ on $N_1\cap \Gamma_1^+$,  
\end{center}
which gives
 \begin{align}\label{eq: 4.6}
\overline{R}_{Q_1}>\overline{S}_{Q_1}> \overline{S}_{R_1}> \overline{R}_{R_1}.
\end{align}
 The last inequality holds because of the fact $\varPi=0$ at $R_1=\partial N_1\cap \Gamma_1^+$.
 
 Now we draw a negative characteristic $\Gamma_1^-$ from $R_1$ to $R_{1*}$ on the boundary $\overline{BA_*}$. We further investigate the analysis of this case by considering the following two subcases. \vspace{0.2 cm}\\
 \textbf{Subcase i :} Assume that $\overline{S}\geq \overline{R}$ holds on every point of $\Gamma_1^-$ then using $\bar{\partial}_- \overline{R}= \overline{R} \varPi$ leads to
$$-\bar{\partial}_- \ln \overline{R}\geq -\Hat{\kappa}.$$
Integrating the above from $R_{1*}$ to $R_1$ we obtain
\begin{equation}\label{eq: 4.7}
\overline{R}_{R_1}\geq \overline{R}_{R_{1*}}e^{-\Hat{\kappa} d}\geq \overline{m}e^{-\Hat{\kappa} d}.
\end{equation}
By using $\eqref{eq: 4.6}$ and $\eqref{eq: 4.7}$ we have
\begin{align}\label{eq: 4.8}
\overline{R}_{Q_1}>\overline{S}_{Q_1}\geq \overline{m}e^{-\Hat{\kappa} d}.
\end{align}
 \textbf{Subcase ii:} Assume that $\overline{S}< \overline{R}$ holds for some point on $\Gamma_1^-$. Then using the fact $\overline{S}>\overline{R}$ at $R_1$ we observe that  $\overline{S}> \overline{R}$ holds for all the points in $N_2\cap \Gamma_1^-$ for some neighbourhood $N_2$ of $R_1$ and $\overline{S}=\overline{R}$ at $Q_2=\partial N_2\cap \Gamma_1^-$.
 
 Then integrating $\eqref{eq: 4.3}$ along negative characteristic from $Q_2$ to $R_1$ and using $\eqref{eq: 4.6}$ yields
 \begin{align}\label{eq: 4.9}
 \overline{R}_{Q_1}>\overline{S}_{Q_1}> \overline{S}_{R_1}> \overline{R}_{R_1}\geq \overline{R}_{Q_2}e^{-\Hat{\kappa}d}=\overline{S}_{Q_2}e^{-\Hat{\kappa}d}.
 \end{align}
 Now from $Q_2$ we draw a positive characteristic $\Gamma_2^+$ up to $Q_2^*$ on the boundary $\overline{BC_*}$. At $Q_2$ we have $\varPi> 0$. Again, if $\varPi \geq 0$ holds on $\Gamma_2^+$ completely then
 $$\overline{S}_{Q_2}>\overline{S}_{Q_2^*}\geq \overline{m},$$
 which in the view of $\eqref{eq: 4.9}$ gives 
 $$\overline{R}_{Q_1}>\overline{S}_{Q_1}\geq \overline{m}e^{-\Hat{\kappa} d}.$$
Otherwise we have a point on $\Gamma_2^+$ such that $\varPi<0$ at that point and
$\bar{\partial}_+S>0$ or $\overline{S}_{Q_2}>\overline{S}_{R_2}$ in the neighbourhood of $Q_2$ where $R_2$ is the point on $\Gamma_2^+$ such that $\varPi>0$ on $\overline{Q_2R_2}\cap \Gamma_2^+$ and $\varPi=0$ at $R_2$. Again using $\eqref{eq: 4.3}$ and $\eqref{eq: 4.9}$ yields
$$\overline{R}_{Q_1}>\overline{S}_{Q_1}\geq \overline{S}_{R_2}e^{-\Hat{\kappa}d}.$$
Since at $R_2$, $\varPi=0$ so we have
$$\overline{R}_{R_2}<\overline{S}_{R_2}.$$
Then we draw a negative characteristic from $R_2$. The repetition of the above process completes the proof of the Lemma.
\end{proof}
\section{Existence of global solution}\label{5}
In this section, we try to extend our ideas to obtain the global solution from the local solution by solving many local Goursat problems in each step of extension. Let us assume that the Goursat problem $\eqref{eq: 2.1}$, $\eqref{eq: 3.2}$ admits a unique $C^1$ solution in $D_{\epsilon}$. Let $\overline{XZ}$ and $\overline{XY}$ be positive($C+$) and negative($C-$) characteristics in $D_{\epsilon}$, respectively. Then, we prescribe 
\begin{align}\label{eq: 5.1}
    (\alpha, \beta, c)=\begin{cases}
    (\alpha|_{XZ}, \beta|_{XZ}, c|_{XZ}), ~~~~~~~~~~~on~~~~~~ \overline{XZ}\\
    (\alpha|_{XY}, \beta|_{XY}, c|_{XY}), ~~~~~~~~~~~on~~~~~~ \overline{XY}
    \end{cases}
\end{align}
where ($\alpha|_{XY}, \beta|_{XY}, c|_{XY}$) and ($\alpha|_{XZ}, \beta|_{XZ}, c|_{XZ}$)  are the values of the solution $(\alpha, \beta, c)$ on $\overline{XY}$ and $\overline{XZ}$, respectively.

We then have the following Lemma.
\begin{l1}\label{l-5.1}
\textit{(Curved quadrilateral building block) If the arc lengths of $\overline{XY}$ and $\overline{XZ}$ are less than $\nu_0$ where $\nu_0=\frac{c_4\sin\left({\frac{\pi}{2}-\epsilon}\right)}{3 M_\epsilon \tan^2\left({\frac{\pi}{4}}+\frac{\epsilon}{2}\right)}\left(\frac{\pi}{4}-\frac{\epsilon}{2}\right)\mu^2(\tau)$,
then the Goursat problem $\eqref{eq: 2.1}$, $\eqref{eq: 5.1}$ admits a global $C^1$ solution on a curved quadrilateral domain bounded by $\overline{XY}, \overline{XZ}, \overline{TY}$ and $\overline{TZ}$ where $\overline{TY}$ is the positive characteristic passing through $Y$, $\overline{TZ}$ is the negative characteristic passing through $Z$. Further, this solution satisfies $\epsilon<\omega<\frac{\pi}{4}+\frac{\epsilon}{2}$.}
\end{l1}
\begin{proof} We know by Lemma $\ref{l-4.1}$ that $\underset{\overline{XZ}}{\sup} \left(\frac{\bar{\partial}_+c}{\sin^2 \omega}\right)\leq M_\epsilon$ and $\underset{\overline{XY}}{\sup}\left( \frac{-\bar{\partial}_-c}{\sin^2 \omega}\right)\leq M_\epsilon$. Then using a similar argument as in the proof of Lemma $\ref{l-4.1}$, one can prove that the solution of Goursat problem $\eqref{eq: 2.1}$, $\eqref{eq: 5.1}$ in the quadrilateral region bounded by $\overline{XY}, \overline{XZ}, \overline{TY}$ and $\overline{TZ}$ satisfies Lemma $\ref{l-3.1}$ and 
\begin{align}\label{eq: 5.2}
    (\bar{\partial}_+c, -\bar{\partial}_-c)\in (0, M_\epsilon]\times (0, M_\epsilon].
\end{align}
Thus using $\eqref{eq: 2.5}$ we can easily obtain 
\begin{align}\label{eq: 5.3}
    0<\bar{\partial}_+\alpha\leq \dfrac{M_\epsilon \tan^2{\omega}}{2 \mu^2(\tau) c_4}, ~~
    -\dfrac{M_\epsilon \tan^2{\omega}}{2 \mu^2(\tau) c_4}\leq -\bar{\partial}_-\beta <0.
\end{align}
We now prove that, when the arc lengths of $\overline{XZ}$ and $\overline{XY}$ are less than $\nu_0$ the solution of the Goursat problem satisfies 
\begin{align}
    \omega< \frac{\pi}{4}+\frac{\epsilon}{2},~~~~~~0<\bar{\partial}_+\alpha\leq \dfrac{M_\epsilon \tan^2({\frac{\pi}{4}+\frac{\epsilon}{2}})}{2 \mu^2(\tau) c_4}, ~~~
    -\dfrac{M_\epsilon \tan^2({\frac{\pi}{4}+\frac{\epsilon}{2}})}{2 \mu^2(\tau) c_4}\leq -\bar{\partial}_-\beta <0.
\end{align} 
We prove it using the method of contradiction. Let us consider an arbitrary point $H$ in the domain of the solution of Goursat problem such that $\omega(H)=\left(\frac{\pi}{4}+\frac{\epsilon}{2}\right)$ and let the positive characteristic passing through $H$ intersects the boundary $\overline{XY}$ at a point $H_1$ and negative characteristic passing through $H$ intersects the boundary $\overline{XZ}$ at a point $H_2$. We draw a straight line passing through $H_1$ having a slope $\tan \alpha(H)$ and another straight line with slope $\tan \beta(H)$ passing through $H_2$. Let us consider that these two lines intersect at a point $J$; see Figure $\ref{fig:my_label5}$. By the estimate $\eqref{eq: 5.3}$ we see that the arcs $\overline{HH_1}$ and $\overline{HH_2}$ lie inside the triangle $JH_1H_2$. Also, their arc lengths are less than the sum of the arc lengths of $\overline{JH_1}, \overline{JH_2}$ and $\overline{H_1H_2}$. Therefore the arc lengths of $HH_1$ and $HH_2$ are less than $\frac{6\nu_0}{\sin({\frac{\pi}{2}-\epsilon})}$. 

Then we integrate $\eqref{eq: 5.3}$ from $H_1$ to $H$ along positive characteristic and from $H$ to $H_2$ along negative characteristic to obtain
\begin{align}\label{eq: 5.4}
    \alpha(H)<\alpha(H_1)+\left(\frac{M_\epsilon \tan^2\left(\frac{\pi}{4}+\frac{\epsilon}{2}\right)}{2 \mu^2(\tau)c_4}\right)\left(\frac{6\nu_0}{\sin({\frac{\pi}{2}-\epsilon})}\right)<\alpha(X)+\frac{\pi}{4}-\frac{\epsilon}{2}\\
     \beta(H)>\beta(H_2)-\left(\frac{M_\epsilon \tan^2\left(\frac{\pi}{4}+\frac{\epsilon}{2}\right)}{2 \mu^2(\tau)c_4}\right)\left(\frac{6\nu_0}{\sin({\frac{\pi}{2}-\epsilon})}\right)>\beta(X)-\frac{\pi}{4}+\frac{\epsilon}{2}
\end{align}
Thus, we have $\omega(H)<\frac{\pi}{4}+\frac{\epsilon}{2}$ which leads to contradiction.

In order to obtain a priori uniform $C^1$ norm estimates of the solution to the Goursat problem we follow the ideas used in the proof of the Lemma $\ref{l-3.1}$ and $\ref{l-3.2}$. Hence using the theory of global classical solutions for quasilinear hyperbolic equations we can obtain the proof of the Lemma $\ref{l-5.1}$ \cite{li1994global}.
\end{proof}
\begin{t1}\label{t-5.1}
\textit{(Global solution)~The semi-hyperbolic patch problem with the boundary data $\eqref{eq: 3.2}$ admits a unique global $C^1$ solution in the region $ABC$ where the curve $\overline{AC}$ is sonic.}
\begin{figure}
\centering
    \includegraphics[width= 6 in]{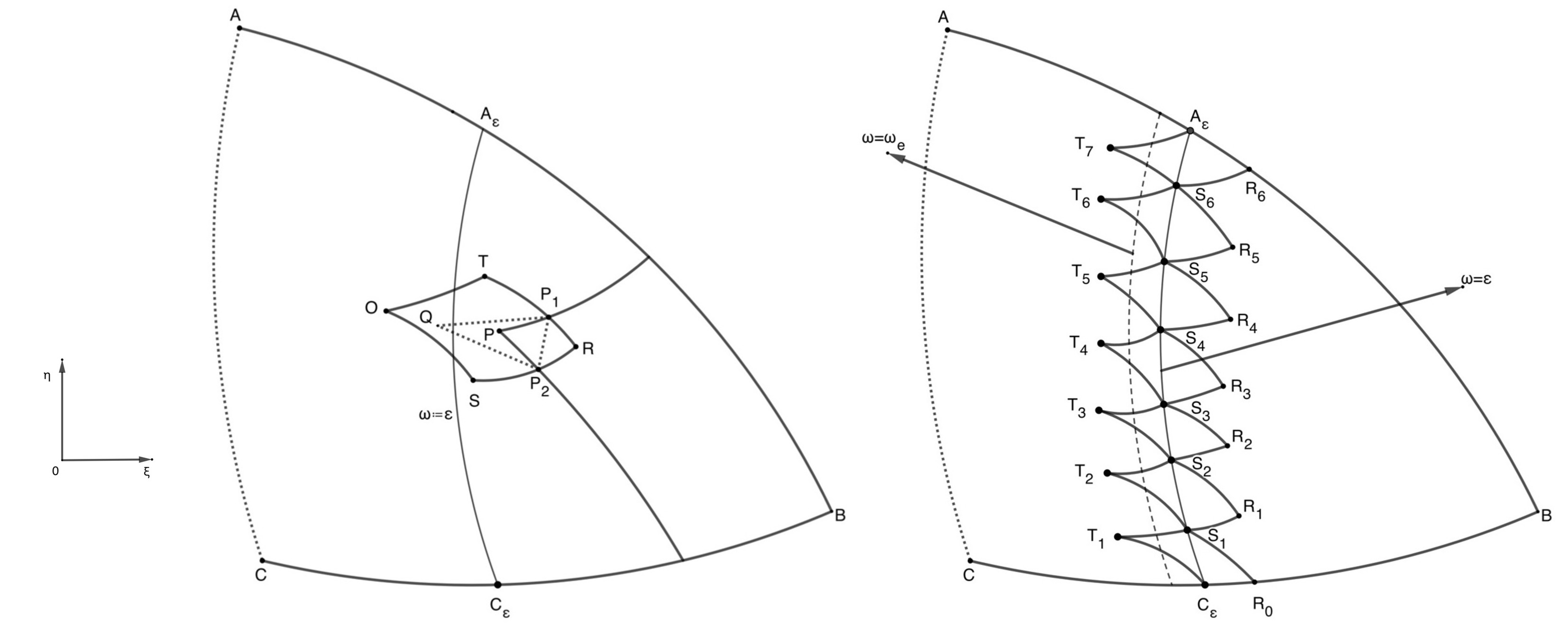}
    \caption{Left: A curved quadrilateral building block; Right: Global solution}
    \label{fig:my_label5}
\end{figure}
\end{t1}
\begin{proof}
Let $Y_0=C_\epsilon, Y_1, Y_2, Y_3,........, Y_n=A_\epsilon$ be $n+1$ different points on the level curve $\omega=\epsilon$. From the point $Y_i$ we draw a $C_+$ characteristic curve which intersects the $C_-$ characteristic curve passing through $Y_{i+1}$ at a point $X_i$, where $i=0,1, 2,....., n-1$. Due to the fact $\bar{\partial}_\pm \omega \neq 0$, the level curve $\omega=\epsilon$ is a non-characteristic curve. Hence, $X_i\neq Y_i$ and $X_i\neq Y_{i+1}$ for any $i=0, 1, 2,......, n-1$. For sufficiently close $Y_i$ and $Y_{i+1}$, the arc lengths of $\overline{X_i Y_i}$ and $\overline{X_i Y_{i+1}}$ are less than $\nu_0$. Therefore, using Lemma $\ref{l-5.1}$, we know that the Goursat problem with the characteristic boundaries $\overline{X_i Y_i}$ and $\overline{X_i Y_{i+1}}$ admits a $C^1$ solution in the quadrilateral domain bounded by $\overline{X_i Y_i}, \overline{X_i Y_{i+1}}, \overline{Y_i Z_{i+1}}$ and $\overline{Y_{i+1} Z_{i+1}}$ where $\overline{Y_i Z_{i+1}}$ is the $C_-$ characteristic curve passing through $Y_i$ and $\overline{Y_{i+1} Z_{i+1}}$ is the $C_+$ characteristic curve passing through $Y_{i+1}$.

Let $Z_0=A$ and $Z_{n+1}=C$, then for every $i=0, 1, 2,.....,n$, there exists a $\omega_i$, $\epsilon<\omega_i<\frac{\pi}{2}$, such that the Goursat problem for system $\eqref{eq: 2.1}$ admits a unique $C^1$ solution in the domain closed by $\overline{Y_i Z_{i+1}}$, $\overline{Y_i Z_{i}}$ and the level curve $\omega(\xi, \eta)=\omega_i$ with $\overline{Y_i Z_i}$ and $\overline{Y_i Z_{i+1}}$ as the characteristic boundaries. Let $\omega_e=\min\{\omega_0, \omega_1, \omega_2,......., \omega_n, \omega(Z_1), \omega(Z_2),......,\omega(Z_n)\}$. Using the fact $\bar{\partial}_+ \omega>0$ we see that $\omega_e>\epsilon$. Then we construct the solution of Goursat problem in the domain $D_{\omega_e}$. Repeating the same process, we can construct the global solution in the whole domain $ABC$ which proves the theorem.
\end{proof}
\section{Shock formation}\label{6}
In this section, we discuss the formation of the envelope for positive characteristics passing through strictly convex curve $\overline{BC}$. Further, we prove that the envelope forms before the sonic points of positive characteristics.
\begin{t1}\label{t-6.1}
\textit{(Envelope formation)~For a given strictly monotonically convex curve $\overline{BC}$, we draw the positive characteristics passing through the curve $\overline{BC}$ which are moving towards downward; see Figure $\ref{fig:my_label6}$. Then positive characteristics form an envelope before their sonic points.}
\begin{figure}
    \centering
    \includegraphics[width= 5 in]{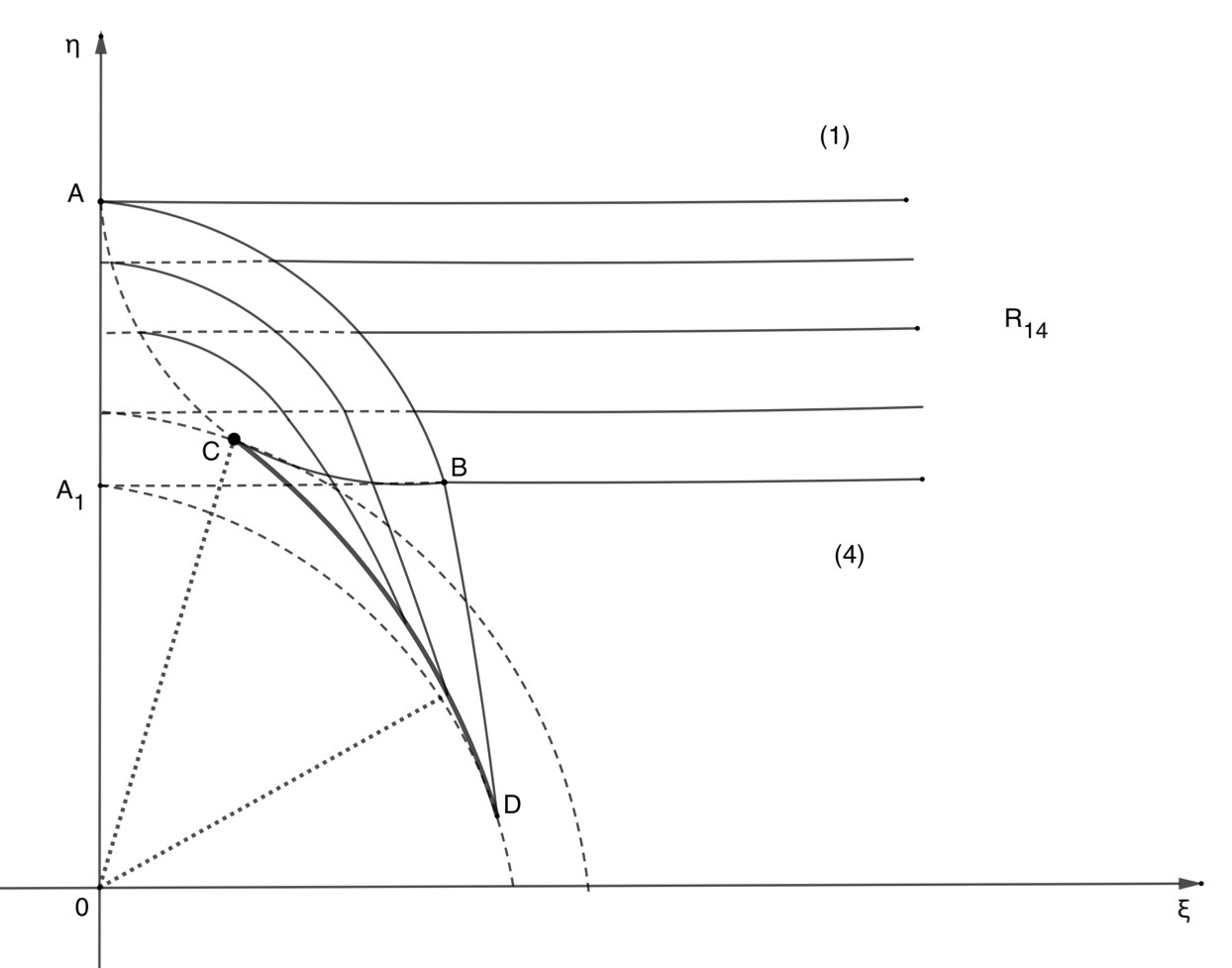}
    \caption{Envelope formation}
    \label{fig:my_label6}
\end{figure}
\end{t1}
\begin{proof}
We exploit $\eqref{eq: 2.6}$ to prove this theorem. Since $\bar{\partial}_+c=0$ in the region of simple waves with the positive characteristics, so from $\eqref{eq: 2.6}$ we obtain
$$\bar{\partial}_+\bar{\partial}_-c=\dfrac{\sin{2 \omega}}{c}\bar{\partial}_-c+\dfrac{(\bar{\partial}_-c)^2}{2\mu^2(\tau)c\cos^2 \omega}.$$
So that
\begin{align}\label{eq: 6.1}
    -\bar{\partial}_+\left(\dfrac{1}{-\bar{\partial}_-c}\right)=\dfrac{\sin{2 \omega}}{c}\left(\dfrac{1}{-\bar{\partial}_-c}\right)-\dfrac{1}{2\mu^2(\tau)c\cos^2 \omega},
\end{align}
which can be written as 
\begin{align}\label{eq: 6.2}
-\bar{\partial}_+\left(\dfrac{1}{-\bar{\partial}_-c} \exp\left({\displaystyle \int_{BD} \dfrac{\sin{2 \omega}}{c}ds}\right)\right)=-\dfrac{1}{2\mu^2(\tau)c\cos^2 \omega}\exp\left({\displaystyle \int_{BD} \dfrac{\sin{2 \omega}}{c}ds}\right).
\end{align}
Using the boundedness of $\sin{2\omega}$ and $c$, we see that $\exp\left({\displaystyle \int_{BD} \frac{\sin{2 \omega}}{c}ds}\right)$ remains bounded on the characteristic extension $\overline{BD}$ which means that the right hand side of $\eqref{eq: 6.2}$ remains negative while on the boundary $\overline{BC}$ we have $\frac{1}{-\bar{\partial}_-c}>0$. Which clearly shows that the function $\frac{1}{-\bar{\partial}_-c}$ is decreasing along the positive characteristics in the direction from $B$ to $D$. Further, the R.H.S of $\eqref{eq: 6.2}$ blows up quadratically at least as shown in \cite{li2011semi}. So $1/(-\bar{\partial}_-c)$ approaches to zero before the characteristic reaches to its sonic point which means that positive characteristics form an envelope before their sonic points.
\end{proof}
 \section{Characteristic decompositions and regularity in partial hodograph plane}\label{7}
 We use the partial hodograph mapping $(\xi, \eta)\longrightarrow (z, t)$ as in \cite{hu2018improved}. We define
 $$t=\cos \omega(\xi, \eta),~~~ z=\phi(\xi, \eta),$$
 where $\phi$ is the potential function used in pseudo-Bernoulli's law $\eqref{eq: 2.2}$.
 
 From the definition of transformation we have the Jacobian $J$ as
 $$J=\dfrac{\partial(z, t)}{\partial(\xi, \eta)}=-\dfrac{(1+\kappa(\tau) \sin^2 
\omega)(\bar{\partial}_+c-\bar{\partial}_-c)}{2 \cos\omega}\neq 0$$ in the entire domain $ABC$.

Let us assume that $A'B'C'$ is the image of the domain $ABC$ in the $z-t$ plane. Then, we transform the normalized derivatives in the new coordinate system $(z, t)$. Using the expression of $\bar{\partial}_\pm$ we obtain
\begin{equation}\label{eq: 7.1}
\begin{aligned}
    \bar{\partial}_+=-\Bigg\{\dfrac{(1+\kappa(\tau)(1-t^2))(1-t^2)\bar{\partial}_+c}{ct}+\dfrac{\sqrt{(1-t^2)^3}}{c}\Bigg\}\dfrac{\partial}{\partial t}+\dfrac{ct}{\sqrt{1-t^2}}\dfrac{\partial}{\partial z},\\
    \bar{\partial}_-=-\Bigg\{\dfrac{(1+\kappa(\tau)(1-t^2))(1-t^2)\bar{\partial}_-c}{ct}+\dfrac{\sqrt{(1-t^2)^3}}{c}\Bigg\}\dfrac{\partial}{\partial t}+\dfrac{ct}{\sqrt{1-t^2}}\dfrac{\partial}{\partial z}.
\end{aligned}
\end{equation}
Using $\eqref{eq: 7.1}$ we obtain 
\begin{equation}\label{eq: 7.2}
    \begin{aligned}
    c_t&=\dfrac{-ct}{(1+\kappa(\tau)(1-t^2))(1-t^2)},~~~~~~~c_z&=\dfrac{t^2-1}{c(1+\kappa(\tau)(1-t^2))}.
    \end{aligned}
\end{equation}
Using the uniform boundedness of $c, ~\kappa(\tau)$ and $t=\cos \omega$, we can verify that $c_t$ and $c_z$ are also uniformly bounded (Recall $c_4<c<c_1$). 

For convenience in the further calculations, we use $R(z, t)=\dfrac{\bar{\partial}_+c}{c}$ and $S(z, t)=\dfrac{\bar{\partial}_-c}{c}$.

By exploiting $\eqref{eq: 7.1}$ in $\eqref{eq: 2.7}$ and $\eqref{eq: 2.8}$ we have 
\begin{equation}\label{eq: 7.3}
    \begin{aligned}
    \begin{cases}
R_t-\dfrac{c f t^2}{S-tg}R_z&=\dfrac{-fR\sqrt{1-t^2}}{S-tg}\Bigg\{\dfrac{2t^2\sqrt{1-t^2}}{c}+\dfrac{R+S}{2\mu^2 t}+\Big\{\tau \kappa'(\tau)-(1+2\kappa(\tau)(1-t^2))\Big\}St \Bigg\},\vspace{0.2 cm} \\ 
S_t-\dfrac{c f t^2}{R-tg}S_z&=\dfrac{-fS\sqrt{1-t^2}}{R-tg}\Bigg\{\dfrac{2t^2\sqrt{1-t^2}}{c}+\dfrac{R+S}{2\mu^2 t}+\Big\{\tau \kappa'(\tau)-(1+2\kappa(\tau)(1-t^2))\Big\}Rt \Bigg\},
\end{cases}
    \end{aligned}
\end{equation}
where $f(t, \tau)=\dfrac{1}{(1+\kappa(\tau)(1-t^2))\sqrt{(1-t^2)^3}}$ and $g(t, \tau)=\dfrac{-(1-t^2)^2 f}{c}$.\vspace{0.2 cm}

We directly compute
\begin{align*}
f_t&=\left(5\kappa(\tau)t\sqrt{(1-t^2)^3}+3t\sqrt{1-t^2}-\kappa'(\tau)(1-t^2)^{\frac{5}{2}} \tau_t\right)f^2,~~~~
f_z=-(1-t^2)^{\frac{5}{2}} f^2 \kappa'(\tau) \tau_z,\\
g_t&=\dfrac{4t(1-t^2)f-(1-t^2)^2f_t}{c}+\dfrac{(1-t^2)^2f}{c^2}c_t,~~~
~~~~~~~~~~~~~~~~~~~~~~~~~g_z=g\left(\dfrac{f_z}{f}-\dfrac{c_z}{c}.\right)
    \end{align*}
where $\tau_t=-\dfrac{\tau \kappa(\tau)}{c}c_t$ and $\tau_z=-\dfrac{\tau \kappa(\tau)}{c}c_z$.

The expressions of $f(t, \tau)$, $g(t, \tau)$, $c_t$, $c_z$, $\tau_t$, $\tau_z$, $f_t$, $f_z$, $g_t$ and $g_z$ clearly show us that these functions are uniformly bounded in the region $A'B'C'$ near sonic boundary, i.e., $t=0$.

Further, we set 
\begin{center}
$\widehat{R}(z, t)=\dfrac{1}{R(z, t)}$ and $\widehat{S}(z, t)=\dfrac{-1}{S(z, t)}$. 
\end{center}

Using Lemma $\ref{l-4.2}$ we see that the functions $\widehat{R}$ and $\widehat{S}$ are uniformly bounded up to sonic curve.\\ 
Now $\eqref{eq: 7.3}$ can be transformed as follows:
\begin{equation}\label{eq: 7.4}
    \begin{aligned}
   \begin{cases}
\widehat{R}_t+\dfrac{c f t^2 \widehat{S}}{1+tg\widehat{S}}\widehat{R}_z=\dfrac{f\sqrt{1-t^2}}{1+tg \widehat{S}}\Bigg\{\dfrac{\widehat{R}-\widehat{S}}{2\mu^2 t}-\dfrac{2t^2\sqrt{1-t^2}\widehat{R}\widehat{S}}{c}+\Big\{\tau \kappa'(\tau)-(1+2\kappa(\tau)(1-t^2))\Big\}\widehat{R}t \Bigg\},\vspace{0.2 cm} \\ 
\widehat{S}_t-\dfrac{c f t^2 \widehat{R}}{1-tg\widehat{R}}\widehat{S}_z=\dfrac{f\sqrt{1-t^2}}{1-tg \widehat{R}}\Bigg\{\dfrac{\widehat{S}-\widehat{R}}{2\mu^2 t}+\dfrac{2t^2\sqrt{1-t^2}\widehat{R}\widehat{S}}{c}+\Big\{\tau \kappa'(\tau)-(1+2\kappa(\tau)(1-t^2))\Big\}\widehat{S}t \Bigg\}.
\end{cases}
    \end{aligned}
\end{equation}
We denote
\begin{align*}
  \Lambda_+=\dfrac{c f t^2 \widehat{S}}{1+tg \widehat{S}},~~~~~  \Lambda_-=\dfrac{-c f t^2 \widehat{R}}{1-tg \widehat{R}}.
    \end{align*}
Further, we denote
$$\partial^\pm=\partial_t+\Lambda_\pm \partial_z,~~~ \Xi=\partial^+ \widehat{R}-\partial^-\widehat{R},~~~ \Theta=\partial^+ \widehat{S}-\partial^-\widehat{S}.$$
Then we have
\begin{align}\label{eq: 7.5}
\Lambda_+ -\Lambda_-=\dfrac{cf(\widehat{R}+\widehat{S})t^2}{(1-tg \widehat{R})(1+tg \widehat{S})},~~~ \widehat{R}_z=\dfrac{\Xi}{\Lambda_+ -\Lambda_-},~~~ \widehat{S}_z=\dfrac{\Theta}{\Lambda_+ -\Lambda_-}.
\end{align}
Now we use the commutator relation \cite{li2006simple}
$$\partial^- \partial^+-\partial^+\partial^-=\dfrac{\partial^- \Lambda_+ -\partial^+ \Lambda_-}{\Lambda_+ -\Lambda_-}(\partial^+-\partial^-),$$
and arrives at
\begin{equation}\label{eq: 7.6}
    \begin{aligned}
\begin{cases}
\partial^+\Xi=\dfrac{\partial^- \Lambda_+ -\partial^+ \Lambda_-}{\Lambda_+ -\Lambda_-}\Xi+(\partial^+\partial^+\widehat{R}-\partial^-\partial^+\widehat{R}), \vspace{0.2 cm}\\
\partial^-\Theta=\dfrac{\partial^- \Lambda_+ -\partial^+ \Lambda_-}{\Lambda_+ -\Lambda_-}\Theta+(\partial^+\partial^-\widehat{S}-\partial^-\partial^-\widehat{S}),    
\end{cases}    
    \end{aligned}
\end{equation}
A straightforward calculation leads to
\begin{equation}\label{eq: 7.7}
    \begin{aligned}
\dfrac{\partial^- \Lambda_+ -\partial^+ \Lambda_-}{\Lambda_+ -\Lambda_-}=\dfrac{2}{t}+h(z, t),    
    \end{aligned}
\end{equation}
where
\begin{align*}
h(z, t)&= \dfrac{f\sqrt{1-t^2}}{(1-tg \widehat{R})(1+tg \widehat{S})}\Bigg\{\dfrac{(\widehat{R}-\widehat{S})g}{2 \mu^2}-\dfrac{2t^3\sqrt{1-t^2}g\widehat{R}\widehat{S}}{c}
+t\Big[\tau \kappa'(\tau)-1-2 \kappa(\tau)(1-t^2)\Big]\Bigg\}\\
&~~~~+\dfrac{(g+tg_t)(\widehat{R}-\widehat{S}+2tg \widehat{R}\widehat{S})}{(1-tg\widehat{R})(1+tg\widehat{S})}+\dfrac{f_t}{f}+tg+\dfrac{cft^3g_z \widehat{R}\widehat{S}}{(1-tg\widehat{R})(1+tg\widehat{S})}.
   \end{align*}
Also, we calculate
\begin{equation}\label{eq: 7.8}
    \begin{aligned}
(\partial^+\widehat{R})_z=f_1\widehat{R}_z+f_2\widehat{S}_z+f_3,   
    \end{aligned}
\end{equation}
where
\begin{align*}
f_1(z, t)&=\dfrac{1}{2t}-\dfrac{(\kappa(\tau)+1)g\widehat{S}-t(1+\kappa(\tau)(2-t^2))(1+tg\widehat{S})}{2(1-t^2)(1+tg\widehat{S})(1+\kappa(\tau)(1-t^2))},\\
&~~~+\dfrac{f\sqrt{1-t^2}}{1+tg \widehat{S}}\Bigg[t(\tau \kappa'(\tau)-1-2\kappa(\tau)(1-t^2))-\dfrac{2t^2\sqrt{1-t^2}\widehat{S}}{c}\Bigg]\\
f_2(z, t)&=\dfrac{-1}{2t}
+\dfrac{(\kappa(\tau)+1)g\widehat{S}-t(1+\kappa(\tau)(2-t^2))(1+tg\widehat{S})}{2(1+\kappa(\tau)(1-t^2))(1-t^2)(1+tg\widehat{S})}
-\dfrac{2t^2f(1-t^2)\widehat{R}}{c(1+tg \widehat{S})}\\
&~~~-\dfrac{fg\sqrt{1-t^2}}{(1+tg\widehat{S})^2}\Bigg[\dfrac{\widehat{R}-\widehat{S}}{2\mu^2}-\dfrac{2t^3\sqrt{1-t^2}\widehat{R}\widehat{S}}{c}+\left(\tau \kappa'(\tau)-\left(1+2\kappa(\tau)(1-t^2)\right)\right)\widehat{R}t^2\Bigg],\\
f_3(z, t)&=\dfrac{(\widehat{R}-\widehat{S})}{2t}\left(\dfrac{\kappa'(\tau)\tau_z}{1+\kappa(\tau)}+\dfrac{f_z}{f(1+tg\widehat{S})}\right)\\&~~~+(\widehat{S}-\widehat{R})\left(\dfrac{\kappa'(\tau)\tau_z}{1+\kappa(\tau)}+\dfrac{f_z}{f(1+tg\widehat{S})}\right)\Bigg[\dfrac{(\kappa(\tau)+1)g\widehat{S}-t(1+\kappa(\tau)(2-t^2))(1+tg\widehat{S})}{2(1+\kappa(\tau)(1-t^2))(1-t^2)(1+tg\widehat{S})}\Bigg]\\
&~~~+\dfrac{f\sqrt{1-t^2}}{1+tg\widehat{S}}\Bigg[-\dfrac{2t^2\sqrt{1-t^2}\widehat{R}\widehat{S}}{\tau \kappa(\tau)c}+\widehat{R}t\Big[\kappa'(\tau)(2t^2-1)+\tau \kappa''(\tau)\Big]\Bigg]\tau_z\\
&~~~+\dfrac{f_z\sqrt{1-t^2}}{(1+tg\widehat{S})^2}\Bigg[-\dfrac{2t^2\sqrt{1-t^2}\widehat{R}\widehat{S}}{c}+\Big[\tau \kappa'(\tau)-\left(1+2\kappa(\tau)(1-t^2)\right)\Big]\widehat{R}t\Bigg]\\
&~~~+\dfrac{fg\widehat{S}\sqrt{1-t^2}}{c(1+tg\widehat{S})^2 }\Bigg[\dfrac{\widehat{R}-\widehat{S}}{2\mu^2}-\dfrac{2t^3\sqrt{1-t^2}\widehat{R}\widehat{S}}{c}+\Big[\tau \kappa'(\tau)-\left(1+2\kappa(\tau)(1-t^2)\right)\Big]\widehat{R}t^2\Bigg]c_z.
\end{align*}
Similarly
\begin{equation}\label{eq: 7.9}
    \begin{aligned}
(\partial^-\widehat{S})_z=g_1\widehat{S}_z+g_2\widehat{R}_z+g_3,
    \end{aligned}
\end{equation}
where
\begin{align*}
 g_1(z, t)&=\dfrac{1}{2t}+\dfrac{(\kappa(\tau)+1)g\widehat{R}+t\left(1+\kappa(\tau)(2-t^2)\right)(1-tg\widehat{R})}{2(1-t^2)(1-tg\widehat{R})\left(1+\kappa(\tau)(1-t^2)\right)}\\
 &+\dfrac{f\sqrt{1-t^2}}{1-tg\widehat{R}}\Bigg[t\left(\tau \kappa'(\tau)-1-2\kappa(\tau)(1-t^2)\right)+\dfrac{2t^2\sqrt{1-t^2}\widehat{R}}{c}\Bigg],\\
 g_2(z, t)&=-\dfrac{1}{2t}-\dfrac{(\kappa(\tau)+1)g\widehat{R}+t\left(1+2\kappa(\tau)(2-t^2)\right)(1-tg\widehat{R})}{2(1-t^2)(1-tg\widehat{R})\left(1+\kappa(\tau)(1-t^2)\right)}+\dfrac{2t^2f(1-t^2)\widehat{S}}{c(1-tg\widehat{R})}\\
 &-\dfrac{fg\sqrt{1-t^2}}{(1-tg\widehat{R})^2}\Bigg[\dfrac{\widehat{R}-\widehat{S}}{2\mu^2}-\dfrac{2t^3\sqrt{1-t^2}\widehat{R}\widehat{S}}{c}-\left(\tau \kappa'(\tau)-\left(1+2\kappa(\tau)(1-t^2)\right)\right)\widehat{S}t^2\Bigg],\\
g_3(z, t)&=\dfrac{\widehat{S}-\widehat{R}}{2t}\left(\dfrac{\kappa'(\tau)\tau_z}{1+\kappa(\tau)}+\dfrac{f_z}{f(1-tg\widehat{R})}\right)\\
&+(\widehat{S}-\widehat{R})\left(\dfrac{\kappa'(\tau)\tau_z}{1+\kappa(\tau)}+\dfrac{f_z}{f(1-tg\widehat{R})}\right)\Bigg[\dfrac{(\kappa(\tau)+1)g\widehat{R}+t\left(1+2\kappa(\tau)(2-t^2)\right)(1-tg\widehat{R})}{2(1-t^2)(1-tg\widehat{R})\left(1+2\kappa(\tau)(1-t^2)\right)}\Bigg]\\
\end{align*}
 \begin{align*}
&+\dfrac{f\sqrt{1-t^2}}{1-tg\widehat{R}}\Bigg[\dfrac{2t^2\sqrt{1-t^2}\widehat{R}\widehat{S}}{\tau \kappa(\tau)c}+\widehat{S}t\Big[\kappa'(\tau)(2t^2-1)+\tau \kappa''(\tau)\Big]\Bigg]\tau_z\\
&+\dfrac{f_z\sqrt{1-t^2}}{(1-tg\widehat{R})^2}\Bigg[\dfrac{2t^2\sqrt{1-t^2}\widehat{R}\widehat{S}}{c}+\left(\tau \kappa'(\tau)-\left(1+2\kappa(\tau)(1-t^2)\right)\right)\widehat{S}t\Bigg]\\
&+\dfrac{fg\widehat{R}\sqrt{1-t^2}}{c(1-tg\widehat{R})^2 }\Bigg[\dfrac{\widehat{R}-\widehat{S}}{2\mu^2}-\dfrac{2t^3\sqrt{1-t^2}\widehat{R}\widehat{S}}{c}-\left(\tau \kappa'(\tau)-\left(1+2\kappa(\tau)(1-t^2)\right)\right)\widehat{S}t^2\Bigg]c_z.
\end{align*}
Usage of $\eqref{eq: 7.5}$ in $\eqref{eq: 7.8}$ and $\eqref{eq: 7.9}$ leads to
\begin{equation}\label{eq: 7.10}
\begin{aligned}
\begin{cases}
\partial^+\partial^+ \widehat{R}-\partial^-\partial^+\widehat{R}=(\Lambda_+-\Lambda_-)(\partial^+\widehat{R})_z=f_1 \Xi+f_2 \Theta+(\Lambda_+-\Lambda_-)f_3,\\
\partial^+\partial^- \widehat{S}-\partial^-\partial^-\widehat{S}=(\Lambda_+-\Lambda_-)(\partial^-\widehat{S})_z=g_1 \Theta+g_2 \Xi+(\Lambda_+-\Lambda_-)g_3.
\end{cases}
\end{aligned}    
\end{equation}
Thus, exploiting $\eqref{eq: 7.6}$, $\eqref{eq: 7.7}$ and $\eqref{eq: 7.10}$, we obtain
\begin{equation}
\begin{aligned}\label{eq: 7.11}
\begin{cases}
\partial^+\Xi=\left(\dfrac{
5}{2t}+\hat{f}_1(z, t)\right)\Xi+\hat{f}_2(z, t)\dfrac{\Theta}{2t}+\hat{f}_3(z, t)\dfrac{t}{2},\vspace{0.2 cm}\\
\partial^-\Theta=\left(\dfrac{
5}{2t}+\hat{g}_1(z, t)\right)\Theta+\hat{g}_2(z, t)\dfrac{\Xi}{2t}+\hat{g}_3(z, t)\dfrac{t}{2},
\end{cases}
\end{aligned}    
\end{equation}
where
\begin{align*}
\hat{f}_1=h+\left(f_1-\dfrac{1}{2t}\right),~~~\hat{f}_2=2tf_2,~~~\hat{f}_3=\dfrac{2}{t}(\Lambda_+-\Lambda_-)f_3,\\
\hat{g}_1=h+\left(g_1-\dfrac{1}{2t}\right),~~~\hat{g}_2=2tg_2,~~~\hat{g}_3=\dfrac{2}{t}(\Lambda_+-\Lambda_-)g_3.
\end{align*}
Noting the expressions of $h, f_i, g_i$ $(i=1, 2, 3)$ and using the Lemma $\ref{l-4.1}$ and $\ref{l-4.2}$, we observe that the functions $\hat{f}_1, \hat{f}_2, \hat{f}_3, \hat{g}_1, \hat{g}_2$ and $\hat{g}_3$ are uniformly bounded near sonic boundary, i.e., $t=0$. Also, we see that $\hat{f}_2\rightarrow -1, \hat{g}_2\rightarrow -1$ as $t\rightarrow 0$. 

Further, we introduce new variables $\widehat{\Xi}=\dfrac{\Xi}{t}, ~\widehat{\Theta}=\dfrac{\Theta}{t}$ to transform the system $\eqref{eq: 7.11}$ into
\begin{align*}
\begin{cases}
\partial^+\widehat{\Xi}=\Bigg\{\dfrac{
3}{2t}+\hat{f}_1(z, t)\Bigg\}\widehat{\Xi}+\hat{f}_2(z, t)\dfrac{\widehat{\Theta}}{2t}+\dfrac{\hat{f}_3(z, t)}{2}, \vspace{0.2 cm}\\
\partial^-\widehat{\Theta}=\Bigg\{\dfrac{
3}{2t}+\hat{g}_1(z, t)\Bigg\}\widehat{\Theta}+\hat{g}_2(z, t)\dfrac{\widehat{\Xi}}{2t}+\dfrac{\hat{g}_3(z, t)}{2},
\end{cases}    
    \end{align*}
or
\begin{equation}\label{eq: 7.12}
    \begin{aligned}
\begin{cases}
\partial^+\left(t^{-\frac{3}{2}}\widehat{\Xi}\right)=t^{-\frac{5}{2}}\Bigg\{t\hat{f}_1(z, t)\widehat{\Xi}+\dfrac{1}{2}\hat{f}_2\widehat{\Theta}+\dfrac{t\hat{f}_3(z, t)}{2}\Bigg\},\vspace{0.2 cm}\\
\partial^-\left(t^{-\frac{3}{2}}\widehat{\Theta}\right)=t^{-\frac{5}{2}}\Bigg\{t\hat{g}_1(z, t)\widehat{\Theta}+\dfrac{1}{2}\hat{g}_2\widehat{\Xi}+\dfrac{t\hat{g}_3(z, t)}{2}\Bigg\}\end{cases}.    
    \end{aligned}
\end{equation}
\subsection{Regularity of functions $R$, $S$ and $W$ in partial hodograph plane}
In this subsection, we are interested to establish the regularity of solution near sonic boundary $\overline{AC}$, i.e., near $t=0$. We use $\eqref{eq: 7.12}$ to derive the regularity results in partial hodograph plane. Let $F=(\phi_1, 0)$ be any point on the line segment $\overline{A'C'}$ where $\overline{A'C'}$ is the image of sonic boundary $\overline{AC}$ in $z-t$ plane. We take a new point $I=(\phi_1, t_m)$ where $t_m$ is very small positive number such that $I=(\phi_1, t_{m})$ remains in the domain $A'B'C'$. Then from the point $I$ we can draw positive and negative characteristic curves $\phi_+(I)$ and $\phi_-(I)$ up to the line segment $\overline{A'C'}$ at $I_1$ and $I_2$, respectively. Since $R, -S$ are uniformly bounded and positive in the domain $ABC$, we see that the functions $\hat{f}_1, \hat{f}_2, \hat{f}_3, \hat{g}_1, \hat{g}_2$ and $\hat{g}_3$ are uniformly bounded near sonic boundary, i.e., $t=0$ in a small subdomain with $\hat{f}_2\rightarrow -1,~\hat{g}_2\rightarrow -1$
as $t\rightarrow 0$. Then for any constant $\varrho \in (0,2]$ we can choose sufficiently small $t_m<1$ such that
\begin{equation}\label{eq: 7.13}
\begin{aligned}
&t|\hat{f_1}|\leq \dfrac{\varrho}{8},~~~~t|\hat{f_3}|\leq \dfrac{\varrho}{4},~~~~|\hat{f_2}|\leq 1+\varrho,\\
&t|\hat{g_1}|\leq \dfrac{\varrho}{8},~~~~t|\hat{g_3}|\leq \dfrac{\varrho}{4},~~~~|\hat{g_2}|\leq 1+\varrho
\end{aligned}
\end{equation}
hold in the domain $II_1I_2$. Let $\Delta(\phi_1, 0)$ be the domain bounded by $II_1$, $II_2$ and the positive and negative characteristics starting from $F=(\phi_1, 0)$. Further, we draw a negative characteristic up to the point $p(z_p, t_p)$ on the boundary $II_1$ and a positive characteristic up to the point $q(z_q, t_q)$ on the boundary $II_2$ starting from an arbitrary point $(z, t)$.
\begin{figure}
    \centering
    \includegraphics[width=4 in]{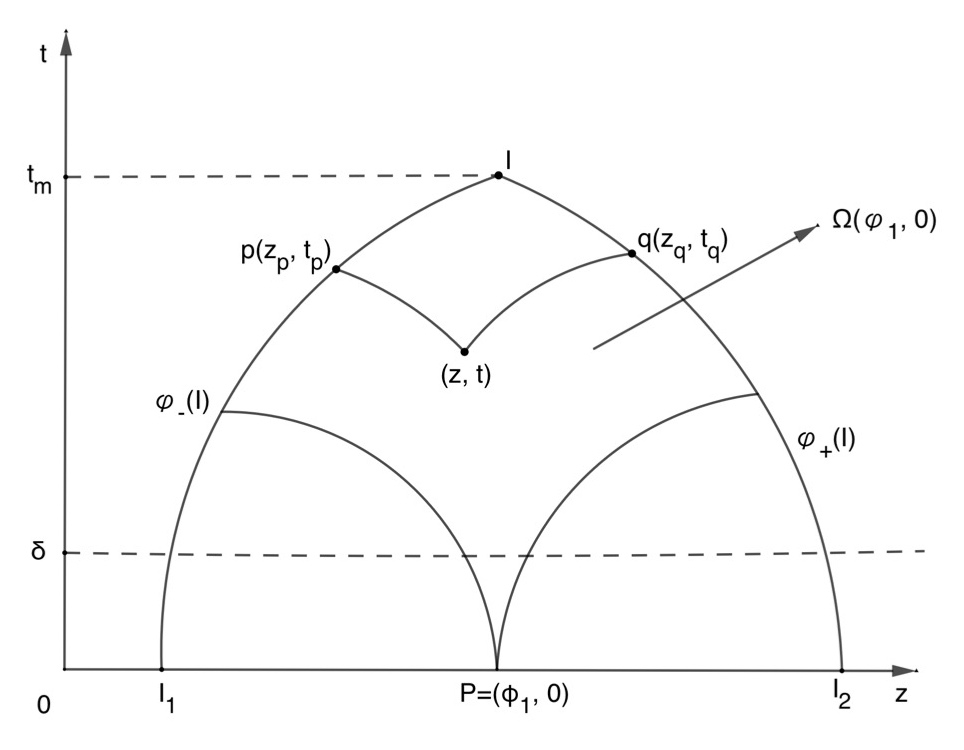}
    \caption{Domain of $\Delta(\phi_1, 0)$.}
    \label{fig:my_label7}
\end{figure}\\
Let us denote
$$M=\max\bigg\{\underset{\Delta(\phi_1, 0)}{\max} \left|\widehat{\Xi}(z_q, t_q)\right|+1,~\underset{\Delta(\phi_1, 0)}{\max}\left|\widehat{\Theta}(z_p, t_p)\right|+1\bigg\}.$$
Since $t_p$ and $t_q$ are strictly positive, so using the definitions of $\widehat{\Xi}$ and $\widehat{\Theta}$ we observe that $M$ is well-defined and uniformly bounded in the domain $\Delta(\phi_1, 0)$.

Now for any fixed $\delta\in (0, t_m)$, define
$$\bar{\Delta}_\delta=\big\{(\phi, t)| \delta \leq t \leq t_m, \phi_-(I)\leq z \leq \phi_+(I)\big\}\cap \Delta(\phi_1, 0);$$
and $M_\delta=\underset{\bar{\Delta}_\delta}{\max}\{|t^\varrho \widehat{\Xi}|, |t^\varrho \widehat{\Theta}|\}$. Then we provide the bounds of $|t^\varrho \widehat{\Xi}|$ and $|t^\varrho \widehat{\Theta}|$ in the following Lemma.
\begin{l1}\label{l-7.1}
\textit{If $(\phi_1, 0)$ is an arbitrary fixed point on the line segment $\overline{A'C'}$  and $\varrho\in (0, 2]$ be any constant. Then there exists a uniform positive constant $\widetilde{M}$ such that the following inequalities hold}
\begin{align*}
\left|t^{\varrho}\widehat{\Xi}\right|\leq \widetilde{M}, ~~~ |t^{\varrho}\widehat{\Theta}|\leq \widetilde{M} ~~ \forall (z, t)\in \Delta(\phi_1, 0).
    \end{align*}
\end{l1}
\begin{proof}~~ Suppose that for every $\delta\in (0, t_m) $, $M_\delta\leq M$ then the Lemma holds true. Otherwise, there exists a $\delta_0\in (0, t_m)$ such that $M_{\delta_0}> M$.

For any point $(z_{\delta_0}, \delta_0)\in \Delta(\phi_1, 0)$, we integrate $\eqref{eq: 7.12}$ along the positive characteristic from $t(\geq \delta_0)$ to $t_q$ and use $\eqref{eq: 7.13}$ to obtain\\
\begin{align*}
&\left|\dfrac{\widehat{\Xi}(z, t)}{t^\frac{3}{2}}\right|=\left|\dfrac{\widehat{\Xi}(z_q, t_q)}{t_q^\frac{3}{2}}+\displaystyle \int_{t}^{t_q}\dfrac{t \hat{f}_1 \widehat{\Xi}+\frac{\hat{f}_2}{2}\widehat{\Theta}+ \frac{t\hat{f}_3}{2}}{t^\frac{5}{2}}dt\right|\\
\leq& \left|\dfrac{\widehat{\Xi}(z_q, t_q)}{t_q^\frac{3}{2}}\right|+\left|\displaystyle \int_{t}^{t_q}\dfrac{\left(\frac{\varrho}{8}\right) \widehat{\Xi}+\left(\frac{1+\varrho}{2}\right)\widehat{\Theta}+\frac{\varrho}{8}}{t^{\frac{5}{2}}}dt\right|\leq \left|\dfrac{\widehat{\Xi}(z_q, t_q)}{t_q^{\frac{3}{2}}}\right|+\dfrac{2+3\varrho}{4}M_{\delta_0}\displaystyle \int_{t}^{t_q}\dfrac{1}{t^{\frac{5}{2}+\varrho}}dt\\
=&\left|\dfrac{\widehat{\Xi}(z_q, t_q)}{t_q^{\frac{3}{2}}}\right|+\dfrac{2+3 \varrho}{4}M_{\delta_0} \dfrac{1}{\left(\frac{3}{2}+\varrho\right)}\left(\dfrac{1}{t^{\frac{3}{2}+\varrho}}-\dfrac{1}{t_q^{\frac{3}{2}+\varrho}}\right)< \dfrac{2+3 \varrho }{6+4\varrho}M_{\delta_0}t^{-\frac{3}{2}-\varrho}
< M_{\delta_0}t^{-\frac{3}{2}-\varrho},
\end{align*}
which proves the inequality
\begin{equation}\label{eq: 7.14}
    \begin{aligned}
 |\widehat{\Xi}(z, t)|_{t={\delta_0}}< M_{\delta_0}{\delta_0}^{-\varrho},
    \end{aligned}
\end{equation}
holds on the line segment $\{t=\delta_0\}\cap \bar{\Delta}_{\delta_0}$. Similarly, we can prove that 
\begin{equation}\label{eq: 7.15}
    \begin{aligned}
 |\widehat{\Theta}(z, t)|_{t={\delta_0}}< M_{\delta_0}{\delta_0}^{-\varrho}.
    \end{aligned}
\end{equation}
 According to $\eqref{eq: 7.14}$ and $\eqref{eq: 7.15}$ we notice that $|t^\varrho \widehat{\Xi}|$ and $|t^\varrho \widehat{\Theta}|$ do not attain the maximum values on the line segment $t=\delta_0$, which means that the maximum value is attained in the interior, i.e., for $\delta_0<t\leq t_m$ in the domain $\bar{\Delta}_{\delta_0}$. The above assertion holds in a larger domain $\bar{\Delta}_{\delta'}$, where $\delta'<\delta_0$. We can extend the domain larger and larger in each step to the whole domain $\Delta(\phi_1, 0)$ to complete the proof of the Lemma.
 \end{proof}
 Next, for any point $(\phi_1, 0)\in \overline{A'C'}$ let $r$ be a positive constant such that $(\phi_1-r, \phi_1+r)\subset A'C'$. Further, assume that the intersection points of negative and positive characteristics passing through points $F_1=(\phi_1-r, 0)$ and $F_2=(\phi_1+r, 0)$ are $F_1^*$ and $F_2^*$, respectively. Then we can extend the inequality in a larger domain $IF_1^*F_1F_2F_2^*$ using the same arguments as in Lemma $\ref{l-7.1}$. Then we have the following Lemma.
 \begin{l1}\label{l-7.2}
 \textit{If $(\phi_1, 0)$ is an arbitrary fixed point on the line segment $\overline{A'C'}$  and $\varrho\in (1, 2]$ be any constant. Then there exists a uniform positive constant $\widetilde{M}$ depending only on $\varrho$ and $r$ such that the following inequalities hold}
\begin{align*}
|t^{\varrho}\widehat{\Xi}|\leq \widetilde{M}, ~~~ |t^{\varrho}\widehat{\Theta}|\leq \widetilde{M} ~~    \forall (z, t)\in IF_1^*F_1F_2F_2^*. 
\end{align*}
\end{l1}
We now prove the uniform boundedness of the function $\widehat{W}=\dfrac{\widehat{R}-\widehat{S}}{t}$.
\begin{l1}\label{l-7.3}
\textit{The function $\widehat{W}$ is uniformly bounded up to the sonic boundary $A'C'$.}
\end{l1}
\begin{proof}~  We exploit the values of $\widehat{R}_t$,  $\widehat{S}_t$, $\widehat{\Xi}$ and $\widehat{\Theta}$ to prove the uniform boundedness of $\widehat{W}$ near sonic boundary, i.e., $t=0$. Therefore, we compute
\begin{equation}\label{eq: 7.16}
    \begin{aligned}
\widehat{W}_t=\theta_1 \widehat{W}+\theta_2,        
    \end{aligned}
\end{equation}
where
\begin{align*}
&\theta_1=\dfrac{(1+\kappa(\tau))g(\widehat{S}-\widehat{R})+2t(1+\kappa(\tau)(2-t^2))+2g(1+\kappa(\tau)(1-t^2))(1-t^2)(\widehat{R}-\widehat{S}+tg \widehat{R}\widehat{S})}{2(1-tg \widehat{R})(1+tg \widehat{S})(1+\kappa(\tau)(1-t^2))(1-t^2)},\\
&\theta_2=\dfrac{f\sqrt{1-t^2}}{(1-tg \widehat{R})(1+tg \widehat{S})}\Bigg[(tg(\Tilde{R}-\Tilde{S})-2)\Bigg\{\dfrac{2t\sqrt{1-t^2}\widehat{R}\widehat{S}}{c}\Bigg\}\Bigg]\\
&~~~~~~+\dfrac{f\sqrt{1-t^2}}{(1-tg \widehat{R})(1+tg \widehat{S})}\Bigg[(\tau \kappa'(\tau)-1-2\kappa(\tau)(1-t^2))(\widehat{R}-\widehat{S}-tg(\widehat{R}^2+\widehat{S}^2))\Bigg]\\
&~~~~~~-\Bigg[\dfrac{\widehat{\Xi}\widehat{S}(1-tg \widehat{R})+\widehat{\Theta}\widehat{R}(1+tg\widehat{S})}{\widehat{R}+\widehat{S}}\Bigg].
\end{align*}
For sufficiently small $t$, we see that $\theta_1$ and $\theta_2$ are uniformly bounded. Thus, integrating $\eqref{eq: 7.16}$ yields that the function $\widehat{W}$ is uniformly bounded near $t=0$.
\end{proof}

Since $\widehat{R}(z, t)=\frac{1}{R(z, t)}$ and $\widehat{S}(z, t)=-\frac{1}{S(z, t)}$ then using the fact that the function $\widehat{W}(z, t)=\frac{\widehat{R}-\widehat{S}}{t}(z, t)$ is uniformly bounded, we can easily obtain the bound of $W(z, t)=\frac{R+S}{t}(z, t)$ which follows that there exists a uniform positive constant $M_1$ such that
\begin{equation}\label{eq: 7.17}
    \begin{aligned}
\left|\dfrac{R+S}{t}(z, t)\right|\leq M_1~~~\forall (z, t)\in A'B'C'.
    \end{aligned}
\end{equation}
 So using $\eqref{eq: 7.3}$ we can prove the uniform boundedness of $|\bar{\partial}^+R|$ and $|\bar{\partial}^-S|$ in the entire domain $A'B'C'$ including the sonic boundary $\overline{A'C'}$, i.e., there exists a uniform positive constant $M_2$ such that
\begin{equation}\label{eq: 7.18}
|\bar{\partial}^+R|<M_2, ~~~~~~~
  |\bar{\partial}^-S|<M_2~~ \forall (z,t)\in A'B'C'.
\end{equation}

We now develop the uniform regularity of functions $R,~S$ and $W$ in partial hodograph plane up to the degenerate line segment $\overline{A'C'}$.
\begin{l1}\label{l-7.4}
\textit{The functions $R(z, t),~S(z, t)$ and $W(z, t)$ are uniformly Lipschitz continuous in the region $A'B'C'$ up to the degenerate line segment $\overline{A'C'}$.}
\end{l1}
\begin{proof}~Using $\eqref{eq: 7.5}$ we obtain
\begin{align}
    R_z=\dfrac{\widehat{\Xi}(R-tg)(S-tg)R^2}{cf(R^2-S^2)}\left(\dfrac{R+S}{t}\right),~~S_z=\dfrac{\widehat{\Theta}(R-tg)(S-tg)S^2}{cf(S^2-R^2)}\left(\dfrac{R+S}{t}\right)
\end{align}

Using $\eqref{eq: 7.17}$, Lemma $\ref{l-4.1}$ and Lemma $\ref{l-4.2}$ we observe that the functions $R_z$ and $S_z$ are uniformly bounded in the domain $A'B'C'$ up to the degenerate line segment $\overline{A'C'}$. Therefore, by $\eqref{eq: 7.18}$ the functions $R_t$ and $S_t$ are also uniformly bounded in the entire region $A'B'C'$.

These observations implies the uniform Lipschitz continuity of the function $R(z,t)$ in the entire domain $A'B'C'$. Similarly, we can prove the uniform Lipschitz continuity of the function $S(z,t)$ in the whole domain $A'B'C'$.

Now we compute $\bar{\partial}^+W$ and $\bar{\partial}^-W$ as follows
\begin{align*}
 \bar{\partial}^+W=-f\sqrt{1-t^2}\Bigg\{\dfrac{2t\sqrt{1-t^2}(R^2+S^2-(R+S)tg)}{c(S-tg)(R-tg)}+\dfrac{RS\left(\tau \kappa'(\tau)-1-2\kappa(\tau)(1-t^2)\right)(R+S-2tg)}{(R-tg)(S-tg)}\Bigg\}\\+\dfrac{cf(S-R)}{(R-tg)(S-tg)}(tS_z)
 +\dfrac{(1+\kappa(\tau)(2-t^2))tW}{(1-t^2)(1+\kappa(\tau)(1-t^2))}
 -\dfrac{(1+\kappa(\tau))(R+S-2tg)(W-g)W}{2(1-t^2)(1+\kappa(\tau)(1-t^2))(S-tg)(R-tg)},
\end{align*}
and 
\begin{align*}
 &\bar{\partial}^-W=-f\sqrt{1-t^2}\Bigg\{\dfrac{2t\sqrt{1-t^2}(R^2+S^2-(R+S)tg)}{c(S-tg)(R-tg)}+\dfrac{RS\left(\tau \kappa'(\tau)-1-2\kappa(\tau)(1-t^2)\right)(R+S-2tg)}{(R-tg)(S-tg)}\Bigg\}\\
 &-\dfrac{cf(S-R)}{(R-tg)(S-tg)}(tR_z)
 +\dfrac{(1+\kappa(\tau)(2-t^2))tW}{(1-t^2)(1+\kappa(\tau)(1-t^2))}
 -\dfrac{(1+\kappa(\tau))(R+S-2tg)(W-g)W}{2(1-t^2)(1+\kappa(\tau)(1-t^2))(S-tg)(R-tg)},
\end{align*}
which clearly indicates that the functions $\bar{\partial}^+W$ and $\bar{\partial}^-W$ are uniformly bounded in the domain $A'B'C'$ up to the degenerate line segment $\overline{A'C'}$. Therefore, a similar proof as of uniform Lipschitz continuity of the functions $R$ and $S$ provide us the uniform Lipschitz continuity of $W$ in $z-t$ plane. Hence the Lemma is proved.
\end{proof}
\section{Regularity of solution in self-similar plane}
\label{8}
Using the results obtained in the preceding section we now prove that the physical variables $(\rho, u, v)(\xi, \eta)$ are uniformly $C^{1, \frac{1}{2}}$ continuous in the region $ABC$ up to the sonic boundary $\overline{AC}$ and also the sonic boundary is $C^{1, \frac{1}{2}}$ continuous. 

To check the regularity of solution in the entire region $ABC$, we first consider the level curves
$$l^\vartheta(\xi, \eta)=1- \sin \omega= \vartheta,$$
where $\vartheta$ is a positive constant. In particular for $\omega= \dfrac{\pi}{2}$, this level curve represents the sonic boundary $\overline{AC}$. 

Now exploiting $\bar{\partial}_+$, $\bar{\partial}_-$ and $\eqref{eq: 2.5}$ we obtain
\begin{equation}\label{eq: 8.1}
\begin{aligned}
l_\xi^\vartheta&=\dfrac{\sin \sigma(1+\kappa(\tau) \sin^2 
\omega)(\bar{\partial}_+c-\bar{\partial}_-c)}{2c}-\dfrac{\sin \omega \cos \sigma (1+\kappa(\tau) \sin^2 \omega)(\bar{\partial}_+c+\bar{\partial}_-c)}{2c \cos \omega}-\dfrac{\sin^2 \omega \cos \sigma }{c}\\
l_\eta^\vartheta&=\dfrac{-\cos \sigma(1+\kappa(\tau) \sin^2 
\omega)(\bar{\partial}_+c-\bar{\partial}_-c)}{2c}-\dfrac{\sin \sigma \sin \omega (1+\kappa(\tau) \sin^2 \omega)(\bar{\partial}_+c+\bar{\partial}_-c)}{2c  \cos \omega}-\dfrac{\sin^2 \omega \sin \sigma}{c}.
\end{aligned}
\end{equation}
From $\eqref{eq: 8.1}$ we obtain
$$(l_\xi^\vartheta)^2+(l_\eta^\vartheta)^2=\left(\dfrac{(1+\kappa(\tau) \sin^2
\omega)(\bar{\partial}_+c-\bar{\partial}_-c)}{2c}\right)^2+\left(\dfrac{\sin \omega (1+\kappa(\tau) \sin^2 \omega)(\bar{\partial}_+c-\bar{\partial}_-c)}{2c \cos \omega}+\dfrac{\sin^2 \omega}{c}\right)^2.$$
Hence, by Lemma $\ref{l-4.1}$ and Lemma $\ref{l-4.2}$ we have
\begin{align}\label{eq: 8.2}
0<\overline{m}^2e^{-2\Hat{\kappa}d}<(l_\xi^\vartheta)^2+(l_\eta^\vartheta)^2\leq (1+\kappa(\tau)M')^2+\Bigg[\dfrac{1}{c_4}+\dfrac{1+\kappa(\tau)}{2}M'\Bigg]^2,
\end{align}
where $M'$ is a uniform positive constant which is the upper bound of $\frac{\bar{\partial}_+c}{c}$ and $-\frac{\bar{\partial}_-c}{c}$.

Now we prove the regularity result in the following four steps.
\subsection{Mapping $(\xi, \eta)\rightarrow (z, t)$ is injective}\label{8.1}
We prove this by the method of contradiction. Let us assume that there exist two distinct points $(\xi_1, \eta_1)$ and $(\xi_2, \eta_2)$ in the region $ABC$ such that $t_1=t_2$ and $z_1=z_2$. 
Which implies that $\cos{\omega(\xi_1, \eta_1)}=\cos{\omega(\xi_2, \eta_2)}$ and $\phi(\xi_1, \eta_1)=\phi(\xi_2, \eta_2)$ such that both the points $(\xi_1, \eta_1)$ and $(\xi_2, \eta_2)$ lie on the same level curve $l^\vartheta=1-\sin{\omega(\xi, \eta)}=\vartheta\geq 0$. Now from $\eqref{eq: 8.1}$ we obtain
\begin{align*}
\nabla \phi.(l^\vartheta_{\eta}, -l^\vartheta_{\xi})=\dfrac{c}{\sin{\omega}}\Bigg\{-\dfrac{(1+\kappa(\tau)\sin^2{\omega}(\bar{\partial}_+c-\bar{\partial}_-c)}{2c}\Bigg\}<0,
\end{align*}
so $\phi(\xi, \eta)$ is strictly monotonically decreasing along each level curve $l^\varpi(\xi, \eta)=\varpi\geq 0$ which contradicts the assumption $\phi(\xi_1, \eta_1)=\phi(\xi_2, \eta_2)$. Hence the mapping is injective.
\subsection{Uniform $C^{\frac{1}{2}}$ continuity of the function $\omega(\xi, \eta)$}
Using $\eqref{eq: 8.1}$ we have
\begin{align}\label{eq: 8.3}
    \cos \omega \omega_\xi=-l^\vartheta_\xi, ~~
    \cos \omega \omega_\eta=-l^\vartheta_\eta
\end{align} 
which gives us 
\begin{align}\label{eq: 8.4}
    |\cos \omega \omega_\xi|+
    |\cos \omega \omega_\eta|\leq 2(1+\kappa(\tau)M')+\Bigg[\dfrac{2}{c_4}+(1+\kappa(\tau))M' \Bigg]\leq C_1
\end{align}
where $C_1$ is a positive constant. From $\eqref{eq: 8.4}$ we obtain
\begin{align*}
    \left|\left(\frac{\pi}{2}-\omega\right)\omega_\xi\right|+
    \left|\left(\frac{\pi}{2}-\omega\right)\omega_\eta\right|\leq \dfrac{\frac{\pi}{2}-\omega}{\sin{(\frac{\pi}{2}-\omega)}}C_1\leq 2C_1
\end{align*}
which follows that the function $(\frac{\pi}{2}-\omega)^2$ is uniformly Lipschitz continuous in $(\xi, \eta)$ plane, which means that for any two points $T_1=(\xi_1, \eta_1)$ and $T_2=(\xi_2, \eta_2)$ in the domain $ABC$ we have
\begin{align}\label{eq: 8.5}
\left|\left(\frac{\pi}{2}-\omega(\xi_2, \eta_2)\right)^2\right|-
    \left|\left(\frac{\pi}{2}-\omega(\xi_1,\eta_1 )\right)^2\right|\leq 2C_1|(\xi_2, \eta_2)-(\xi_1, \eta_1)|
\end{align}
Since $\frac{\pi}{4}< \omega < \frac{\pi}{2}$, we observe that $\left(\dfrac{\pi}{2}-\omega(\xi_2, \eta_2)\right)$ and $\left(\dfrac{\pi}{2}-\omega(\xi_1, \eta_1)\right)$ are positive so that
\begin{align*}
&~~~~\left|\left(\frac{\pi}{2}-\omega(\xi_2, \eta_2)\right)^2-\left(\frac{\pi}{2}-\omega(\xi_1, \eta_1)\right)^2\right|\\
    &=\left|\left(\frac{\pi}{2}-\omega(\xi_2, \eta_2)\right)-
    \left(\frac{\pi}{2}-\omega(\xi_1, \eta_1)\right)\right|.\left|\left(\frac{\pi}{2}-\omega(\xi_2, \eta_2)\right)+
    \left(\frac{\pi}{2}-\omega(\xi_1, \eta_1)\right)\right|\\
    &\geq \left|\left(\frac{\pi}{2}-\omega(\xi_2, \eta_2)\right)-
    \left(\frac{\pi}{2}-\omega(\xi_1, \eta_1)\right)\right|.\left|\left(\frac{\pi}{2}-\omega(\xi_2, \eta_2)\right)-
    \left(\frac{\pi}{2}-\omega(\xi_1, \eta_1)\right)\right|\\
    &=|\omega(\xi_2, \eta_2)-\omega(\xi_1, \eta_1)|^2
\end{align*}
 Hence by $\eqref{eq: 8.5}$, we obtain
 \begin{align}\label{eq: 8.6}
|\omega(\xi_2, \eta_2)-\omega(\xi_1, \eta_1)|\leq \sqrt{2C_1}|(\xi_2, \eta_2)-(\xi_1, \eta_1)|^\frac{1}{2}     
 \end{align}
which proves the uniform $C^{\frac{1}{2}}$ continuity of the function $\omega(\xi, \eta)$ in the entire domain $ABC$.

\subsection{Uniform regularity of functions $\frac{\bar{\partial}_+c}{c}, \frac{\bar{\partial}_-c}{c}$ and $\left(\frac{\bar{\partial}_+c+\bar{\partial}_-c}{c \cos{\omega}}\right)$ in the entire region $ABC$ of $(\xi, \eta)$ plane}
From subsection $\ref{8.1}$ we know that the mapping $(\xi, \eta)\rightarrow (z, t)$ is injective, so for any two distinct points $(\xi_1, \eta_1)$ and $(\xi_2, \eta_2)$ we have two different images, say $T_1'=(z_1, t_1)$ and $T_2'=(z_2, t_2)$ in the region $A'B'C'$. 
Then in this subsection we prove the uniform $C^{\frac{1}{2}}$ continuity of the functions $\widetilde{R}(\xi, \eta):=R(z, t),~\widetilde{S}(\xi, \eta):=S(z, t)$ and $\widetilde{W}(\xi, \eta):=W(z, t)$ where $\widetilde{R}(\xi, \eta)=\bar{\partial}_+c/c$, $\widetilde{S}(\xi, \eta)=\bar{\partial}_-c/c$ and $\widetilde{W}(\xi, \eta)=(\bar{\partial}_+c+\bar{\partial}_-c)/c \cos\omega$ in $\xi-\eta$ plane.

We use the uniform Lipschitz continuity of $R(z,t)$ from Lemma $\ref{l-7.4}$ to obtain
\begin{align}\label{eq: 8.7}
   |\widetilde{R}(\xi_2, \eta_2)-\widetilde{R}(\xi_1, \eta_1)| &=|R(z_2, t_2)-R(z_1, t_1)|\leq C_2|(z_2, t_2)-(z_1, t_1)|\nonumber \\
   &=C_2\Big\{\left(\cos{\omega}(\xi_2, \eta_2)-\cos{\omega}(\xi_1, \eta_1)\right)^2+\left(\phi(\xi_2, \eta_2)-\phi(\xi_1, \eta_1)\right)^2\Big\}^\frac{1}{2}
\end{align} 
for some uniform positive constant $C_2$.

Now using $\eqref{eq: 8.6}$ we have
\begin{align*}
    |\cos{\omega}(\xi_2, \eta_2)-\cos{\omega}(\xi_1, \eta_1)|
    &\leq \left|2\sin\left({\dfrac{\omega(\xi_2, \eta_2)-\omega(\xi_1, \eta_1)}{2}}\right)\right|\\
    &\leq |\omega(\xi_2, \eta_2)-\omega(\xi_1, \eta_1)| \leq \sqrt{2C_1}|(\xi_2, \eta_2)-(\xi_1, \eta_1)|^\frac{1}{2}.   
\end{align*}
So, the uniform Lipschitz continuity of potential function $\phi(\xi, \eta)$ in the entire region $ABC$ and $\eqref{eq: 8.7}$ yields
\begin{align*}
   |\widetilde{R}(\xi_2, \eta_2)-\widetilde{R}(\xi_1, \eta_1)|&\leq C_2\Big\{2C_1|(\xi_2, \eta_2)-(\xi_1, \eta_1)|+C_3^2|(\xi_2, \eta_2)-(\xi_1, \eta_1)|^2\Big\}^\frac{1}{2}\\
   &\leq C_4|(\xi_2, \eta_2)-(\xi_1, \eta_1)|^\frac{1}{2}  
\end{align*}
where $C_3$ and $C_4$ are positive constants such that $C_4= \max\{C_2\sqrt{2C_1}, C_2C_3\}$. 

The above result concludes the uniform $C^\frac{1}{2}$ continuity of the function $\widetilde{R}(\xi, \eta)$ in the entire region $ABC$. In the same manner, we can prove that the functions $\widetilde{S}(\xi, \eta)$ and $\widetilde{W}(\xi, \eta)$ are also uniformly $C^{\frac{1}{2}}$ continuous.
\subsection{Uniform regularity of the solution $(\rho, u, v)$ in the entire region $ABC$ of $(\xi, \eta)$ plane}
Now we derive the uniform regularity of solution $(\rho, u, v)$ and sonic boundary $\overline{AC}$ using the uniform regularity of functions $\widetilde{R}, ~\widetilde{S}$ and $\widetilde{W}$ in the entire region $ABC$ in $(\xi, \eta)$ plane. To prove this we first prove the uniform regularity of $\sigma(\xi, \eta)$.

Using $\eqref{eq: 2.5}$ we obtain
\begin{align}\label{eq: 8.8}
    \begin{cases}
\sigma_\xi=\dfrac{\sin\sigma \sin \omega}{c}+\dfrac{\left(\sin \sigma \cos \omega(\widetilde{R}+\widetilde{S})-\cos \sigma \sin \omega(\widetilde{R}-\widetilde{S})\right)\kappa(\tau)}{2},    \vspace{0.1 cm}\\
\sigma_\eta=-\dfrac{\cos\sigma \sin \omega}{c}-\dfrac{\left(\cos \sigma \cos \omega(\widetilde{R}+\widetilde{S})+\sin \sigma \sin \omega(\widetilde{R}-\widetilde{S})\right)\kappa(\tau)}{2}    
    \end{cases}
\end{align}
which clearly shows that the functions $\sigma(\xi, \eta), \cos{\sigma(\xi, \eta)}$ and $\sin{\sigma(\xi, \eta)}$ are uniformly Lipschitz continuous in the whole region $ABC$. Also, using the fact that $\widetilde{R}, \widetilde{S}, \omega\in C^{\frac{1}{2}}$ and $\eqref{eq: 8.8}$ we see that $\sigma(\xi, \eta)\in C^{1, \frac{1}{2}}$ and eventually $\cos{\sigma(\xi, \eta)}, ~\sin{\sigma(\xi, \eta)} \in C^{1, \frac{1}{2}}$. Using this and $\eqref{eq: 8.1}$ one can prove that the function $\sin \omega(\xi, \eta)$ is uniformly $C^{1,\frac{1}{2}}$ continuous in the whole domain $ABC$. 

Further, using the pseudo-Bernoulli's law we notice that $\rho(\xi, \eta)$ is uniformly $C^{1, \frac{1}{2}}$ continuous in the whole region $ABC$ which means that the function $c(\xi, \eta)$ is uniformly $C^{1, \frac{1}{2}}$ continuous. The expressions for pseudo-velocities lead to uniform $C^{1, \frac{1}{2}}$ continuity of the functions $u(\xi, \eta)$ and $v(\xi, \eta)$ in the entire region $ABC$.

Also from $\eqref{eq: 8.1}$ and $\eqref{eq: 8.2}$ we notice that $l^\vartheta_\xi$ and $l^\vartheta_\eta$ are $C^{\frac{1}{2}}$ continuous and $(l^\vartheta_\xi)^2+(l^\vartheta_\eta)^2$ is bounded. Therefore, the sonic boundary $\overline{AC}$ is $C^{1, \frac{1}{2}}$ continuous. 

We summarize the regularity results in the following theorem.
\begin{t1}\label{t-8.1}
~\textit{If the angle $\beta_C\in(-\frac{\pi}{2}, 0)$, then the Goursat problem $\eqref{eq: 2.1}$ and $\eqref{eq: 3.2}$ admits a global smooth solution in the region $ABC$ where the curve $\overline{AC}$ is the sonic boundary. Further, this solution is uniformly $C^{1, \frac{1}{2}}$ continuous up to the sonic boundary $\overline{AC}$ and the sonic boundary $\overline{AC}$ is $C^{1, \frac{1}{2}}$ continuous.}
\end{t1}
\section{Conclusions}\label{9} 
In this paper, we considered a special case of convex pressure and proved that the global solution to the semi-hyperbolic patch problem for two-dimensional compressible Euler equations with van der Waals gas exists and these solutions are uniformly $C^{1, \frac{1}{2}}$ continuous and also the sonic boundary is $C^{1, \frac{1}{2}}$ continuous. The study of semi-hyperbolic patch problem opens door to extend the solution into the subsonic domain, which we will try to tackle in future.\\\\

\nocite{*}
\bibliographystyle{elsarticle-num}  \bibliography{name}
\end{document}